# VALIDITY OF HEAVY TRAFFIC STEADY-STATE APPROXIMATIONS IN GENERALIZED JACKSON NETWORKS

By David Gamarnik and Assaf Zeevi[1]

*MIT and Columbia University*

We consider a single class open queueing network, also known as a generalized Jackson network (GJN). A classical result in heavy-traffic theory asserts that the sequence of normalized queue length processes of the GJN converge weakly to a reflected Brownian motion (RBM) in the orthant, as the traffic intensity approaches unity. However, barring simple instances, it is still not known whether the stationary distribution of RBM provides a valid approximation for the steady-state of the original network. In this paper we resolve this open problem by proving that the re-scaled stationary distribution of the GJN converges to the stationary distribution of the RBM, thus validating a so-called "interchange-of-limits" for this class of networks. Our method of proof involves a combination of Lyapunov function techniques, strong approximations and tail probability bounds that yield tightness of the sequence of stationary distributions of the GJN.

**1. Introduction.** Jackson networks are one of the most fundamental objects in the theory of stochastic processing networks, owing to the remarkable simplicity of the product form stationary distribution of the queue lengths. Product form results have since been established for wider classes of open and closed queueing networks, however, these results typically hinge on exponential distributions being imposed on the interarrival and service times. Relaxing these assumptions to i.i.d. interarrival and service times with general distributions significantly complicates the analysis and, barring special cases, does not lead to closed form expressions for the steady-state queue lengths distribution. From an application standpoint, however, such assumptions typically constitute a more adequate model of "real world" systems. Consequently, the study of open single class queueing networks with

Received November 2004; revised July 2005.
[1]Supported in part by NSF Grant DMI-04-97652.

Bernoulli type routing and i.i.d. interarrival and service times, hereafter referred to as generalized Jackson networks (GJN), has been the focus of vigorous academic research.

In the absence of closed form solutions for the steady-state queue length distribution in GJN, efforts have mostly centered on the development of various bounds and approximations. Most notably, much of the recent work on single and multiclass queueing networks has focused on *fluid models* for stability analysis, *diffusion approximations* for performance analysis and the relationship between the two methods and objectives. Below we provide a brief overview of these developments in the context of GJN.

In terms of diffusion approximations, Reiman's seminal paper [33] proved that, as the traffic intensity approaches one, the normalized vector of queue length processes converges in distribution to a Brownian motion which is reflected at the boundaries of the positive orthant. The normalization in Reiman's heavy-traffic limit theorem corresponds to the functional central limit theorem, that is, time is accelerated linearly and space is compressed by a square-root factor. The limit process, known as reflected Brownian motion (RBM), was first constructed in earlier work of Harrison and Reiman [18]. The main appeal of the RBM process is that it provides a parsimonious (two moment) mesoscopic approximation of the underlying discrete dynamics in the GJN. Reiman's theorem provides a rigorous justification for the use of RBM as an approximation to the original (appropriately normalized) queue length process in the GJN.

To approximate the steady-state behavior of the network, it is first necessary to derive conditions under which both the GJN and RBM admit a stationary distribution. To this end, Sigman [35] and Down and Meyn [15] were the first to prove that the GJN possesses a unique stationary distribution if and only if the traffic intensity at every station is less than unity. A similar result was later proved by Dai [12] using a remarkable connection between the stochastic and fluid models of queueing networks; see also [37]. (Fluid limits of GJN were first characterized by Johnson [22] and Chen and Mandelbaum [9, 10].) Despite the fact that the stability conditions in GJN have such a simple form, exact computation of the stationary distribution is not possible except in the case of product form networks.

It turns out that the stability condition for RBM is quite similar to that of the underlying GJN. In particular, it was shown by Harrison and Williams [19] that a (unique) stationary distribution exists for the RBM, roughly speaking, if and only if the "traffic intensity" at each station is less than unity (a precise definition will be given in the following section). For further details and references, the reader is referred to the survey paper by Williams [38]. Although the stationary distribution of an RBM cannot be typically computed in closed form, with the exception of special cases involving a skew-symmetry condition (see, e.g., [19] and further discussion in Section 4.4), it can be solved for numerically. This computation typically exploits the partial differential equation, also known as the basic adjoint relation, that characterizes the stationary distribution. The algorithms of Dai and Harrison [13] and Shen et al. [36] build on this idea and describe different approaches to compute the stationary distribution by numerically solving the basic adjoint partial differential equation.

Given the apparent solidarity between GJN and its associated RBM approximation, as well as the fact that steady-state computations can be more easily performed for the latter, it has become commonplace to approximate the steady-state of GJN with that of the associated RBM. It is important to note that this approximation *is not* rigorously justified and, in particular, does not follow from Reiman's [33] heavy-traffic process level approximation. To date, this correspondence has only been established in some specific problem instances, starting with Kingman's seminal work [27, 28] on the $GI/GI/1$ queue. Specifically, Kingman proved that the distribution of the diffusion-scaled stationary waiting time converges to an exponential distribution, which is identical to the steady-state of the RBM obtained as a process limit from the normalized queue length. (For extensions of this result that cover dependent input processes, see, e.g., [34].) In the context of a single station multiclass queueing network with feedback, the workload process, after diffusion scaling, also converges to a one-dimensional RBM. The steady-state mean value of this RBM is known to be equal to the expected workload in the original model under diffusion scaling. Thus, in this example the validity of the RBM steady-state approximation is established thanks to the fact that both of the above expectations can be computed in closed form. Harrison [17] proved that, in a series of two single-server queues, the scaled steady-state vector of queue lengths converges weakly to a random vector with the steady-state distribution of the associated RBM. Kaspi and Mandelbaum [23] proved a similar convergence for a *closed* GJN. (This case is quite distinct from the open network setting since the state space in the former is compact and, hence, establishing tightness of the sequence of re-scaled queue lengths is straightforward.) A very interesting result in the context of large deviations theory was established by Majewski [31]. He showed that, in a feedforward-type multiclass queueing network with deterministic service times, the large deviations rates corresponding to the network converge to the large deviations rates of the associated RBM when the network is in heavy traffic. While these results are suggestive of the fact that the stationary distribution of RBM may indeed provide a rigorous approximation to that of the underlying GJN, the recent book of Chen and Yao [11] mentions that this remains an open problem. Resolving this question would constitute an important stepping stone on route to proving the validity of RBM steady-state approximations in other classes of queueing networks.

4  D. GAMARNIK AND A. ZEEVI

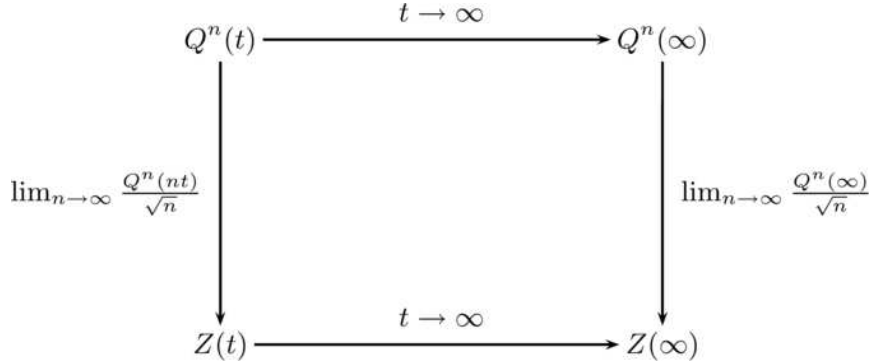

FIG. 1. *The interchange-of-limits diagram establishing the validity of heavy-traffic steady-state approximations.*

The main goal of the present paper is to establish this conjectured correspondence between the steady-state of GJN and its Brownian approximation. In particular, we prove that the stationary distribution of the vector of queue lengths in the GJN, under diffusion scaling, converges to the stationary distribution of the associated RBM as the traffic intensity approaches unity (see Theorem 8). This result establishes the validity of an "interchange of limits" argument. Specifically, by first taking heavy-traffic limits, one arrives at an RBM approximation from which, by letting time go to infinity, one obtains the steady-state approximation to the GJN. What we show is that this order of limits can be interchanged. In particular, we prove that the same stationary distribution is obtained by first passing to steady-state in the GJN, and subsequently applying diffusion scaling and letting the system utilization approach one. This argument is represented symbolically in the diagram given in Figure 1, where $Q^n(t)$ represents the queue length process in heavy-traffic, and $Z(t)$ represents the associated RBM. Using the above, we also prove that the normalized moments with respect to the stationary distribution of the GJN converge to the respective steady-state moments of the RBM (see Corollary 2). Finally, we establish the above interchange argument for the sojourn times (see Theorem 10), thus proving a steady-state version of Reiman's [33] "snapshot principle."

The main technical subproblem that we solve in order to establish our main result is that of tightness of the normalized queue length vector in steady-state. Roughly speaking, we show that, as the traffic intensity $\rho_j \to 1$ in every station $j = 1, \ldots, J$, the family of scaled queue length vectors $((1-\rho_j)Q_j(t))_j$ is tight (see Theorem 7). This establishes the existence of a probability measure which is a weak limit of this family. It is then a matter of simple argument to show that this measure must be the unique stationary distribution of the associated RBM. Proving tightness of the family of normalized steady-state queue lengths hinges on the following ingredients: (i) existence of a pathwise solution to the Skorohod problem (the oblique reflection mapping); (ii) existence of a "fluid-scale" linear Lyapunov function for the workload; and (iii) probability bounds for strong approximations of the primitive processes, that is, arrivals, services and routing.

On route to our main result, we derive exponential bounds on the stationary distribution of general Markov processes admitting a suitable Lyapunov function (see Theorems 5 and 6). While such bounds are excessive for purposes of proving tightness of the stationary distribution under diffusion scaling, they are potentially quite useful for performance analysis of GJN. (For recent work that derives tail bounds in the heavy-tailed context for GJN, and more general classes of monotone separable networks, see [1].) These bounds and their derivation are in the spirit of the results obtained by Meyn and Tweedie in several papers and summarized in their book [32]. Unlike the aforementioned results of Meyn and Tweedie [32] (see also the queueing application in [14]), our work does not use an infinitesimal drift characterization of the Markov process in question. Rather, we derive our bounds on the basis of the behavior of the Markov process on "large" time scales. This approach is more intuitive, supports a more direct application of the theory to the GJN context, and circumvents several technical issues. In terms of strong approximations, we rely on similar ideas to those used by Horváth [21], however, our main focus is on establishing uniform integrability of diffusion-scaled processes. The blend of methods we employ might be of independent interest and, in a parallel ongoing work, we use these techniques to derive explicit upper and lower bounds on steady-state queue lengths in GJN under both heavy-tailed and light-tailed assumptions on the primitives.

The main technical assumption we impose on the primitives of the GJN is the following. Both the distribution of interarrival and service times are assumed to have a finite conditional exponential moment generating function (a precise definition is given in what follows). We note that "standard" heavy-traffic limit theory, including Reiman's theorem, typically only assumes the existence of a second moment (more precisely $2 + \delta$ moments) for the interarrival and service time distributions. In that sense our conditions are clearly much stronger than necessary. Nevertheless, these conditions facilitate the mathematical derivations and, moreover, lead to novel exponential tail bounds on the steady-state of the GJN (see Theorem 6). Extending our results to the more general case assuming only existence of the first two moments is an important direction for future work. Note that the RBM corresponding to the heavy-traffic limit of the underlying GJN has a steady-state with exponential tails; see [19]. This suggests that the existence of exponential moments for the GJN primitives should result in a more accurate RBM-based steady-state approximation. To this end, it is worth noting that, even in the simple case of the $GI/GI/1$ queue with heavy-tailed service times, the RBM exponential steady-state is known to provide a poor approximation of the Pareto-like stationary distribution of the original queue length process. In Section 5 we discuss briefly conditions under which we believe our analysis should extend to the case of GJN with primitives having Pareto-like distributions.

The rest of the paper is organized as follows. Section 2 provides background and preliminaries on GJN, the Skorohod mapping and heavy-traffic approximations. Section 3 derives bounds on the stationary distribution of general state-space Markov chains using Lyapunov arguments. Section 4 details the main results. Specifically, the main theorem is stated in Section 4.2. Concluding remarks are given in Section 5, and all proofs are collected in Appendix.

**2. The queueing network model and its heavy-traffic approximation.** The set up and notation follow closely that in [11], Chapter 7. We use $I$ to denote a $J \times J$ identity matrix, and $e$ the $J$-dimensional vector of ones. (All the vectors are assumed to be column vectors, unless otherwise stated.) Put $e^i$ to be the $i$th unit vector in $\mathbb{R}^J$. The transposed of a matrix $P$ (vector $e$) will be denoted $P'$ ($e'$). For every $z \in \mathbb{R}^J$, the norm $\|z\|$ corresponds to the $L_1$ metric: $\|z\| = \sum_{1 \leq j \leq J} |z_j|$. Given a random variable $X$ and a probability measure $\nu$, $X \sim \nu$ means $X$ is distributed according to $\nu$. For a sequence of probability measures $\nu_n, n = 1, 2, \ldots, \nu_n \Rightarrow \nu$ is used to denote weak convergence of $\nu_n$ to a limit probability measure $\nu$. A similar notation is used to denote weak convergence of a sequence of random variables or processes.

2.1. *Generalized Jackson networks*: *model description and probabilistic assumptions.*

*Description of the network and probabilistic assumptions.* The network consists of $J$ stations denoted for simplicity by $1, 2, \ldots, J$. Jobs arrive from external sources to each station according to independent renewal processes with interarrival times given by $a_j(0), a_j(1), a_j(2), \ldots$ for $j = 1, \ldots, J$. The interarrival times $a_j(1), a_j(2), \ldots$ are assumed to have a common distribution $F_{A,j}$ for $j = 1, \ldots, J$. The first interarrival time $a_j(0)$ is assumed to have the same distribution function, only conditioned on the event $a_j(0) \geq a_{0,j}$, where $a_{0,j} \geq 0$ measures the time elapsed since the last arrival prior to time $t = 0$. This is a given deterministic or random value and is considered to be part of the data. Thus, the distribution of the arrival process is completely specified by a pair $(F_{A,j}, a_{0,j})$. Put $a_0 := (a_{0,1}, \ldots, a_{0,J})$. Also, let $\alpha_j = 1/\mathbb{E}[a_j(1)]$ be the external arrival rate into station $j$, and put $\alpha = (\alpha_1, \ldots, \alpha_J)$ to be the vector of external arrival rates. It is not excluded that $\alpha_j = 0$ for some stations $j$, meaning there is no external arrival into these stations.

We denote by $\mathcal{J} = \{j : \alpha_j > 0\}$. Let $A_j(t) = \sup\{n : \sum_{0 \leq l \leq n} a_j(l) \leq t\}$ denote the corresponding counting renewal process. For every $t_2 \geq t_1 \geq 0$, we put $A_j(t_1, t_2) := A_j(t_2) - A_j(t_1)$. The coefficient of variation $\mathrm{var}[a_j^2(1)]\alpha_j^2$ is denoted by $c_{a,j}^2$, for $j \in \mathcal{J}$.

Denote by $v_j(0), v_j(1), \ldots$ the sequence of service times received by jobs in station $j$ and $V_j(n) = \sum_{0 \leq l \leq n} v_j(l)$. The service times $v_j(1), v_j(2), \ldots$, are assumed to be i.i.d. with a distribution function $F_{S,j}$. As in the case of interarrival times, the first service time $v_j(0)$ is assumed to follow the same distribution, only conditioned on $v_j(0) \geq v_{0,j}$. Here $v_{0,j} \geq 0$ is a (possibly random) value which stands for the cumulative service time elapsed for the job currently in service in station $j$ at time $t = 0$. This value is considered a part of the data. For consistency, we assume that $v_{0,j} > 0$ only if the queue length in station $j$ at time 0 is positive. We let $S_j(t) = \sup\{n : V_j(n) \leq t\}$ denote the associated renewal process and let $m_j = \mathbb{E}[v_j(1)]$ denote the mean service time in station $j$. Then, $\mu_j = 1/m_j$ is the service rate in station $j$, and $\mu = (\mu_1, \ldots, \mu_J)$ is the vector of service rates. We let $M$ be the diagonal matrix with $m = (m_1, \ldots, m_J)$ as diagonal entries. The coefficient of variation $\mathrm{var}[v_j^2(1)]\mu_j^2$ is denoted by $c_{s,j}^2$, for $j = 1, \ldots, J$.

We assume that both the interarrival and service times have a finite moment generating function in some neighborhood of the origin. Moreover, to handle the residual service and interarrival times, we assume the following condition holds: there exists $\theta^* > 0$ such that, for every $j$,

(1) $$\sup_{z \in \mathbb{R}_+} \mathbb{E}[\exp(\theta^*(a_j(1) - z))|a_j(1) > z] < \infty,$$

(2) $$\sup_{z \in \mathbb{R}_+} \mathbb{E}[\exp(\theta^*(v_j(1) - z))|v_j(1) > z] < \infty,$$

where $j \in \mathcal{J}$ in (1) and $j = 1, 2, \ldots, J$ in (2). In particular, $\mathbb{E}[\exp(\theta^* a_j(1))]$ and $\mathbb{E}[\exp(\theta^* v_j(1))]$ are finite as well.

The routing decisions at each station are assumed to be of Bernoulli type and parameterized by a sub-stochastic $J \times J$ matrix $P = (p_{jk})_{1 \leq j,k \leq J}$. The entries of $P$ satisfy $p_{ij} \geq 0$, for every $i$, $\sum_j p_{ij} \leq 1$ and the spectral radius of $P$ is strictly less than unity. For each station $j$, let $\psi^j = \{\psi^j(1), \psi^j(2), \ldots\}$ be an i.i.d. sequence of routing decisions with common distribution $P_j$, where $P_j$ is the $j$th column of $P$. Specifically, $\mathbb{P}(\psi^j(l) = e^k) = p_{jk}$ and $\mathbb{P}(\psi^j(l) = 0) = 1 - \sum_{1 \leq k \leq J} p_{jk}$, the latter indicates the fact that the completed job leaves the network. Let $R^j(0) = 0$ and $R^j(n) = \sum_{1 \leq l \leq n} \psi^j(l)$. Then the $k$th component, $R_k^j(n)$, of $R^j(n)$ is the total number of jobs (out of $n$) which, upon completing their service at station $j$, are routed to station $k$.

Let $\lambda = (\lambda_1, \ldots, \lambda_J)$ denote the unique solution to the traffic equation $\lambda = \alpha + P'\lambda$. That is, $\lambda = [I - P']^{-1}\alpha$. Then $\rho = M^{-1}\lambda = M^{-1}[I - P']\alpha =$ $(\rho_j)_{1 \leq j \leq J}$ is the vector of traffic intensities. We let $\rho^* = \max_j \rho_j$ denote the bottleneck traffic intensity. We will primarily consider the case when

$$\rho^* < 1. \tag{3}$$

This is also referred to as the stability condition for GJN, a terminology which is justified by Theorem 2 below.

*Systems dynamics.* The number of jobs which are either in service or waiting in station $j$ at time $t$ is denoted by $Q_j(t)$; adhering with standard terminology, we refer to $Q_j = (Q_j(t) : t \geq 0)$ as the *queue length process* at station $j$. Put $Q = (Q_1, \ldots, Q_J)$. Let $B_j(t)$ denote the cumulative time that station $j$ was busy processing work during the interval $[0, t]$. In particular, $0 \leq B_j(t) \leq t, B(0) = 0$. Put $B_j = (B_j(t) : t \geq 0)$ and let $B = (B_1, \ldots, B_J)$; we refer to this as the *busy time process*. Let $I_j(t) = t - B_j(t)$ denote the cumulative idle time at station $j$ in the interval $[0, t]$. Function $I_j = (I_j(t) : t \geq 0)$ describes the *idle time process* in station $j$. At every station $j$ the jobs are assumed to be processed in First-In-First-Out (FIFO) fashion.

The following equations that describe the system dynamics tie together the processes introduced above:

$$Q_j(t) = Q_j(0) + A_j(t) + \sum_{1 \leq k \leq J} R_j^k(S_k(B_k(t))) - S_j(B_j(t)) \tag{4}$$

and

$$I_j(t) = \int_0^t \mathbb{1}\{Q_j(s) = 0\} \, ds, \tag{5}$$

for all $t \geq 0, j = 1, 2, \ldots, J$. We consider the centered process $X = (X_j(t) : t \geq 0)$ with

$$
\begin{aligned}
X_j(t) = {} & Q_j(0) + \left(\alpha_j + \sum_{1 \leq i \leq J} \mu_i p_{ij} - \mu_j\right) t + (A_j(t) - \alpha_j t) \\
& + \sum_{1 \leq i \leq J} p_{ij}(S_i(B_i(t)) - \mu_i B_i(t)) - (S_j(B_j(t)) - \mu_j B_j(t)) \\
& + \sum_{1 \leq i \leq J} (R_j^i(S_i(B_i(t))) - p_{ij} S_i(B_i(t))),
\end{aligned}
\tag{6}
$$

$j = 1, 2, \ldots, J$, and the process $Y_j = (Y_j(t) : t \geq 0)$, where

$$Y_j(t) = \mu_j I_j(t) = \mu_j(t - B_j(t)). \tag{7}$$

The process $X_j$ is often referred to as the "free process" which summarizes the net input of jobs to station $j$. The process $Y_j$ summarizes system idleness in normalized units. Alternatively, this process "regulates" the queue lengths from taking negative values by increasing whenever the queue length process is zero. The dynamics in (4)–(5) can now be represented in the following succinct form:

$$Q(t) = X(t) + [I - P']Y(t) \geq 0, \tag{8}$$

$$Y(0) = 0 \text{ and } Y(t) \text{ is nondecreasing in } t, \tag{9}$$

$$\int_0^t Q_j(\tau)\,dY_j(\tau) = 0 \qquad \text{for all } t \geq 0, 1 \leq j \leq J. \tag{10}$$

The equations (8)–(10) describe the system dynamics $(Q, Y)$ as the solution to the *Skorohod reflection problem* associated with the GJN. A formal definition of the Skorohod problem is given in Section 2.2 below.

*A Markovian representation.* We now provide a Markov process representation for the underlying queue length processes. The queue length process $Q$ is clearly not Markovian because of the residual interarrival and service times. To create a Markovian structure, we augment the state descriptor as follows. For each time $t \geq 0$ and $j \in \mathcal{J}$, we let $\hat{a}_j(t)$ denote the time elapsed since the most recent arrival to station $j$ that occurred prior to time $t$. That is, $\hat{a}_j(t) := t - \sum_{n \leq A_j(t)} a_j(n)$, and we put $\hat{a}(t) = (\hat{a}_1(t), \ldots, \hat{a}_J(t))$, where, by definition, $\hat{a}_j(t) = 0$ for $j \notin \mathcal{J}$ for all $t$. Similarly, $\hat{v}_j(t)$ denotes the time elapsed since the commencement of the most recent service prior to time $t$ at station $j$, or if the server is idle at time $t$, this variable is equal to zero. Put $\hat{v}(t) = (\hat{v}_1(t), \ldots, \hat{v}_J(t))$. The *extended* process $\bar{Q} = (\bar{Q}(t) : t \geq 0)$ with $\bar{Q}(t) = (Q(t), \hat{a}(t), \hat{v}(t))$ is a Markov process. That is, for every $t_1 > t \geq 0$ and every $(z, a, v) \in \mathbb{Z}_+^J \times \mathbb{R}_+^{2J}$, the distribution of $\bar{Q}(t_1)$ conditioned on $\bar{Q}(t) = (z, a, v)$ is independent of $\bar{Q}(t')$ for all $t' < t$. (In fact, the process $\bar{Q}$ is known to be a strong Markov process; see [11] for further details.) By definition, $\bar{Q}(0) = (Q(0), \hat{a}(0), \hat{v}(0)) = (Q(0), a_0, v_0)$. For brevity, we denote the extended state space of the process $\bar{Q}$ by $\mathcal{X} := \mathbb{Z}_+^J \times \mathbb{R}_+^{2J}$. The probability distribution of the strong Markov process $\bar{Q}(t)$ for every $t \geq 0$ is completely specified by $(Q(0), a_0, v_0, F_{A,j}, F_{S,j}, P)$. We denote this 6-tuple by $\Xi$ and call it the parameter vector for the corresponding generalized Jackson network. For brevity, we often refer to the underlying GJN simply as $\Xi$. A probability distribution $\pi$ on $\Xi$ is defined to be *stationary* if $\bar{Q}(0) \sim \pi$ implies that $\bar{Q}(t) \sim \pi$ for all $t \geq 0$. [The stationary distribution for the Markov process $(\bar{Q}(t) : t \geq 0)$ is then constructed in the usual manner from $\pi$ and the transition kernel.] When $\bar{Q}(t) \sim \pi$ we say that the GJN is in *steady state* with stationary distribution having $\pi$ as its marginal; this distribution is not necessarily unique. For future purposes, we simply use $\pi$ when making reference to the stationary distribution of the Markov process in question. We will also say that the network is *stable* if there exists at least one stationary distribution $\pi$.

2.2. *The Skorohod mapping, fluid model and stability of GJN.* This subsection reviews some results concerning the Skorohod problem and fluid analysis of GJN that will be useful in what follows. The presentation is based on Chapter 7 in [11].

*Skorohod problem and the oblique reflection mapping.* Let $\mathcal{D}[0,\infty)$ denote the space of right-continuous functions with left limits taking values in $\mathbb{R}^J$. Since all limit processes considered in this paper have continuous sample paths, it is not necessary to introduce the Skorohod topology, and it suffices to consider the space $\mathcal{D}$ endowed with the uniform topology; see [6] for further discussion. We now cite a general result, originally derived in [18] in the context of reflected Brownian motion, that supports the representation of GJN dynamics given in (8)–(10); the version below is adapted from [11], Theorem 7.2.

THEOREM 1 (Oblique reflection mapping). *For every function $x \in \mathcal{D}[0,\infty)$ such that $x(0) \in \mathbb{R}_+^J$, there exists a unique pair of nonnegative functions $(y, q) = (\Psi_1(x), \Psi_2(x)) \in \mathcal{D}[0,\infty)$ such that $y_j(t)$ is a nondecreasing function for every $j = 1, 2, \ldots, J$, $y_j(0) = 0$, and $(x, y, q)$ jointly satisfy*

$$q(t) = x(t) + [I - P']y(t), \tag{11}$$

$$\int_0^t q_j(\tau) \, dy_j(\tau) = 0 \qquad \text{for all } t \geq 0 \text{ and } j = 1, 2, \ldots, J. \tag{12}$$

*Moreover, the mappings $\Psi_1, \Psi_2 : \mathcal{D}[0,\infty) \to \mathcal{D}[0,\infty)$ are Lipschitz continuous, in particular, for every $x_1, x_2 \in \mathcal{D}[0,\infty)$ and every $t \geq 0$,*

$$\sup_{0 \leq \tau \leq t} \|q_1(\tau) - q_2(\tau)\| \leq R_1 \sup_{0 \leq \tau \leq t} \|x_1(\tau) - x_2(\tau)\| \tag{13}$$

*and*

$$\sup_{0 \leq \tau \leq t} \|y_1(\tau) - y_2(\tau)\| \leq R_2 \sup_{0 \leq \tau \leq t} \|x_1(\tau) - x_2(\tau)\|, \tag{14}$$

*where $R_1, R_2 > 0$ depend only on the matrix $P$.*

REMARK 1. The processes $q$ and $y$ are referred to as the *reflected process* and the *regulator process*, respectively. For further details and explicit identification of the Lipschitz constants $R_1, R_2$, see, for example, the recent book by Whitt [39], Chapter 14.2.3.

REMARK 2. The existence and uniqueness of the pathwise solution to the Skorohod problem hinges on the matrix $P$ being substochastic. This, in turn, implies that the *reflection matrix* $I - P'$ is Minkowskii, that is, has all positive diagonal elements, nonpositive off-diagonal elements, and inverse in which all entries are nonnegative.

*The fluid model and stochastic stability.* Given any $z \in \mathbb{R}_+^J$, let $x_z(t) = z + (\alpha - (I - P')\mu)t$ and let $(y, q) = (\Psi_1(x), \Psi_2(x))$. The triplet $(x, y, q)$ is called a *fluid model* of the GJN. Let $w = e'[I - P']^{-1}$. Note that each coordinate $w_1, \ldots, w_J$, of the vector

$$w = e'[I - P']^{-1}$$
$$= e'[I + P' + (P')^2 + (P')^3 + \cdots]$$

is greater than one. Consider the process

(15) $$(w'q(t) : t \geq 0).$$

It is well known that $w'q(t)$ is Lipschitz continuous (see [11]). The process $(w'q(t) : t \geq 0)$ can be interpreted, essentially, as the total amount of "fluid work" present in the system at time $t$. An important property of the fluid model is the amount of time needed to drain existing work from the system when all exogenous inflows are "turned off" at time $t$. For the fluid model of a GJN network, we have the following result, which can be found, for example, in [11].

PROPOSITION 1 (Fluid stability of GJN). *The vector valued function $q(t)$ is differentiable everywhere on $\mathbb{R}_+$ except for a set of points which has a Lebeasgue measure zero. Moreover, $w'\dot{q}(t) \leq -\min_{1 \leq j \leq J} \mu_j(1 - \rho_j)$ whenever $q(t) \neq 0$ and $q(t)$ is differentiable. As a result, if condition* (3) *holds, then for any initial fluid level $q(0) = z$, there exists*

(16) $$\tau \leq \frac{w'z}{\min_{1 \leq j \leq J} \mu_j(1 - \rho_j)}$$

*such that $q(t) = 0$ for all $t \geq \tau$.*

The function $w'q(t)$ serves as a Lyapunov function for the GJN network. In particular, the result articulated above asserts that the fluid model is stable, as it drains in finite time from any initial fluid level. Given the results of Dai [12] and Stolyar [37] linking stochastic and fluid stability, Proposition 1 can be used to establish stability of the underlying GJN. To this end, Sigman [35], Down and Meyn [15], Dai [12] and Dai and Meyn [14] prove the following theorem.

THEOREM 2 (Stochastic stability of GJN). *Suppose condition* (3) *holds for the GJN $\Xi$. Then, there exists a stationary distribution $\pi$ for the Markovian process $\bar{Q}(t)$. Under additional technical regularity (e.g., Assumptions A1–A3 in* [14]*), the stationary distribution is unique and for every starting state $\bar{Q}(0) = (z, a, v) \in \mathcal{X}$, the distribution of $\bar{Q}(t)$ converges to $\pi$ as $t \to \infty$.*

The Lyapunov function $w'q(\cdot)$ and the fluid model of the GJN will play an important role in deriving our main results in what follows.

2.3. *Heavy-traffic limits and reflected Brownian motion.*

*The GJN in heavy-traffic.* We now introduce the generalized Jackson network in the parameter regime corresponding to heavy-traffic assumptions. We consider a sequence of networks with traffic intensities $\rho_j^n, n = 1, 2, \ldots$, each approaching unity at the same rate as $n \to \infty$. To arrange for that, we first consider a network such that the traffic intensity $\rho_j = 1$ in every station $j$, and denote this network by $\Xi$. (At this point the reader should recall that $\Xi$ encodes the 6-tuple given in Section 2.1 that parameterizes the GJN.) In particular, $\lambda = [I - P']^{-1}\alpha = \mu$. Fix a vector of strictly positive constants $\kappa^0 = (\kappa_1^0, \ldots, \kappa_j^0)$. Then, we consider a sequence of networks $\{\Xi^n\}$ in which: (i) the service times $v_j$ (as well as the residual service times) and routing decisions $\psi^j$ have the same distribution as in the original network $\Xi$; (ii) the interarrival times for the $n$th network in the sequence are defined to be $a_j^n = a_j(1 - \kappa_j^0/\sqrt{n})$, and the residual interarrival times remains unchanged. Hence, the vector of arrival rates is given by $\alpha^n = (\alpha_1^n, \ldots, \alpha_J^n)$. As a result, the vector of effective arrival rates is $\lambda^n = [I - P']^{-1}\alpha^n$, which we write as $[I - P']^{-1}\alpha - n^{-1/2}[I - P']^{-1}\kappa^0 = \lambda - n^{-1/2}\kappa$, where $\kappa := [I - P']^{-1}\kappa^0 = (\kappa_1, \ldots, \kappa_J)$. Then

$$\text{(17)} \qquad \rho_j^n = (1 - \kappa_j/\sqrt{n}) \qquad \text{for } j = 1, \ldots, J.$$

Equivalently, $\sqrt{n}(1 - \rho_j^n) = \kappa_j$ for $j = 1, \ldots, J$. We also obtain for the bottleneck traffic intensity $\rho^{*n} = (1 - \kappa^*/\sqrt{n})$, where $\kappa^* = \min \kappa_j$. The GJN in heavy-traffic refers to the sequence of generalized Jackson networks $\{\Xi^n\}$, and the corresponding sequences of processes are denoted $A^n(t)$, $B^n(t)$, $I^n(t)$, $X^n(t)$, $Q^n(t)$, $Y^n(t)$. Note that the processes $S(t)$ and the routing processes do not depend on $n$, so we use the notation $S(t), R^j(l)$ for these elements in the network $\Xi^n$. The sequence of stationary distributions corresponding to $\Xi^n$ will be denoted by $\{\pi^n\}$ henceforth. These distributions exist by Theorem 2 due to the fact that condition (3), $\rho^{*n} = \max_j \rho_j^n < 1$, holds for all $n$. When the stationary distributions are not unique, we denote by $\{\pi^n\}$ an arbitrary sequence of stationary distributions corresponding to $\{\Xi^n\}$.

*Reflected Brownian motion in the orthant.* The importance of the heavy-traffic regime stems from the ability to approximate the underlying queue lengths and workload processes by a reflected Brownian motion (RBM) which we define below.

DEFINITION 1 (Reflected Brownian motion in the orthant). Let $z \in \mathbb{R}^J$, $P$ be a substochastic $J \times J$ matrix, and let $W = (W(t) : t \geq 0)$ be a $J$-dimensional Brownian motion with initial state $z$, drift vector $\beta \in \mathbb{R}^J$ and covariance matrix $\Gamma \in \mathbb{R}^{J \times J}$. Then, the reflected Brownian motion (RBM) with parameters $(\beta, \Gamma, I - P')$ and initial state $z$ is the unique solution to the Skorohod problem (11), (12) with input process $W$. In particular, the RBM $Z = (Z(t): t \geq 0)$ is given by $Z = \Psi_1(W)$, and the nondecreasing regulator process $Y = (Y(t): t \geq 0)$ is given by $Y = \Psi_2(W)$. Thus, the pair $(Z, Y)$ satisfies

$$(18) \qquad Z(t) = W(t) + [I - P']Y(t),$$

$$(19) \qquad \int_0^t Z_j(s)\, dY_j(s) = 0 \qquad \text{for all } t \geq 0 \text{ and } j = 1, 2, \ldots, J.$$

The RBM process $Z$ has the following behavior: in the interior of the nonnegative orthant $\mathbb{R}_+^J$ evolves like a Brownian motion with drift $\beta$ and diffusion matrix $\Gamma$; and when the process hits the $j$th face, $F_j = \{x \in \mathbb{R}_+^J : x_j = 0\}$, it is reflected instantaneously toward the interior of the orthant in the direction determined by the $j$th column of the matrix $I - P'$. Roughly speaking, the $j$th component of the regulator process $Y$ measures the cumulative "effort" in maintaining the Brownian motion confined to the nonnegative orthant.

The stationary distribution of the RBM (when it exists) is one of the primary focuses of the present paper. To that end, the following result, established in [19], will be important in what follows.

THEOREM 3 (RBM stationary distribution). *The RBM $(\beta, \Gamma, I - P')$ possesses a stationary distribution $\pi_{\mathrm{RBM}}$ if and only if $[I - P']^{-1}\beta < 0$. When this distribution exists, it is unique. Moreover, for every initial distribution, the marginal distribution of RBM converges weakly to $\pi_{\mathrm{RBM}}$ as time diverges to infinity.*

*Heavy-traffic scaling and diffusion limits.* For the sequence of GJN, $\{\Xi^n\}$, defined earlier, put

$$(20) \qquad \begin{aligned} \hat{X}^n(t) &= X^n(nt)/\sqrt{n}, \\ \hat{Q}^n(t) &= Q^n(nt)/\sqrt{n}. \end{aligned}$$

The above correspond to a *diffusion scaling* of the free process $X$ and the queue length process $Q$. In particular, the normalization in (20) corresponds to accelerating time and compressing space in accordance with the functional central limit theorem (FCLT). The following theorem, due to Reiman [33], is one of the most fundamental results in the theory of diffusion approximations of queueing networks. The theorem was originally proven under the assumption that $Q^n(0) = 0$, however, the more general result stated below also holds true (see, e.g., [11]).

THEOREM 4 (Process limit weak convergence of GJN to RBM). *Consider a sequence of GJN, $\{\Xi^n\}$, with $\rho^n$ satisfying (17), and assume that $Q^n(0)/\sqrt{n} \Rightarrow Z(0)$, where $Z(0)$ is distributed according to some probability measure $\pi_0$ on $\mathbb{R}_+^J$. Then, for every $t > 0$, the sequence of processes $(Q^n(nt')/\sqrt{n} : 0 \leq t' \leq t)$ converges weakly to the RBM $(Z(t') : 0 \leq t' \leq t)$ with respect to the topology of uniform convergence in $\mathcal{D}[0,t]$. The RBM has initial distribution $Z(0) \sim \pi_0$, and it has parameters $(\beta, \Gamma, I - P')$, where the drift vector is given by $\beta = -(I - P')M^{-1}\kappa$, and the entries of the covariance matrix $\Gamma$ are given by*

$$\Gamma_{k\ell} = \sum_{j=1}^{J} \mu_j [p_{jk}(\delta_{k\ell} - p_{j\ell}) + c_{s,j}^2 (p_{jk} - \delta_{jk})(p_{j\ell} - \delta_{j\ell})] + \alpha_k c_{a,k}^2 \delta_{k\ell},$$

*where $\delta_{jk} = 1$ if $j = k$ and is zero otherwise.*

REMARK 3. The limiting covariance structure can be inferred from (48) in Section 5; for further details, see [11], Theorem 7.29.

REMARK 4. Note that the expression for the RBM drift vector implies that $[I - P']^{-1}\beta = -M^{-1}\kappa < 0$, which is the condition for the existence and uniqueness of the stationary distribution of RBM, and corresponds to the condition $\rho_n^* < 1$ in $\Xi^n$.

REMARK 5. An important point concerning the RBM approximation is that it can give rise to closed form solutions for its stationary distribution, while the underlying GJN may be intractable in that regard. Harrison and Williams [19, 20] provide necessary and sufficient conditions under which the stationary distribution of RBM is separable with exponentially distributed marginals. An example that illustrates this is given in Section 4.

## 3. Lyapunov functions and steady-state bounds for Markov processes.

3.1. *Background and definitions.* Let $\Xi = (\Xi(t) : t \geq 0)$ be a continuous time Markov process defined on a complete metrizable state space $\mathcal{X}$, which, for our purposes, is Euclidean space, equipped with Borel $\sigma$-algebra $\mathcal{B}$. For the special case of generalized Jackson network, $\Xi(t) = \bar{Q}(t) = (Q(t), \hat{a}(t), \hat{v}(t))$.

Let $\mathbb{P}_x$ be the probability distribution under which $\mathbb{P}(\Xi(0) = x) = 1$, and put $\mathbb{E}_x[\cdot] := \mathbb{E}[\cdot | \Xi(0) = x]$, that is, the expectation operator w.r.t. the probability distribution $\mathbb{P}_x$. Let $\mathbb{P}_\pi$ denote the probability distribution under which $\Xi(0)$ is distributed according to $\pi$, and put $\mathbb{E}_\pi[\cdot]$ to be the expectation operator w.r.t. this distribution. A probability distribution $\pi$ defined on $\mathcal{B}$ is said to be a *stationary distribution* if for every bounded continuous function $f : \mathcal{X} \to \mathbb{R}$,

$$\mathbb{E}_\pi[f(\Xi(0))] = \mathbb{E}_\pi[f(\Xi(t))], \tag{21}$$

for all $t \geq 0$. The following definition plays an important role in our analysis.

DEFINITION 2 (Lyapunov function). A function $\Phi : \mathcal{X} \to \mathbb{R}_+$ is said to be a Lyapunov function with drift size parameter $-\gamma < 0$, drift time parameter $t_0 > 0$ and exception parameter $K$, if

$$\sup_{x \in \mathcal{X} : \Phi(x) > K} \{\mathbb{E}_x \Phi(\Xi(t_0)) - \Phi(x)\} \leq -\gamma. \qquad (22)$$

A function $\Phi : \mathcal{X} \to \mathbb{R}_+$ is defined to be a geometric Lyapunov function with a geometric drift size $0 < \gamma < 1$, drift time $t_0$ and exception parameter $K$, if

$$\sup_{x \in \mathcal{X} : \Phi(x) > K} \{(\Phi(x))^{-1} \mathbb{E}_x \Phi(\Xi(t_0))\} \leq \gamma. \qquad (23)$$

Given any function $\Phi : \mathcal{X} \to \mathbb{R}_+$ and $t \geq 0$, let

$$\phi(t) = \sup_{x \in \mathcal{X}} \Phi^{-1}(x) \mathbb{E}_x \Phi(\Xi(t)). \qquad (24)$$

It is not excluded that $\phi(t) = \infty$, however, the interesting case of course is when $\phi(t)$ is finite. Given $\theta > 0$, let

$$L_1(\theta, t) = \sup_{x \in \mathcal{X}} \mathbb{E}_x [\exp(\theta(\Phi(\Xi(t)) - \Phi(x)))] \qquad (25)$$

and

$$L_2(\theta, t) = \sup_{x \in \mathcal{X}} \mathbb{E}_x [(\Phi(\Xi(t)) - \Phi(x))^2 \exp(\theta(\Phi(\Xi(t)) - \Phi(x))^+)] \qquad (26)$$

for $t \geq 0$.

3.2. *Bounds on the stationary distribution.* Lyapunov functions provide a useful tool for obtaining bounds on stationary distributions of Markov chains and processes. In this subsection we are set to do just that. The bounds we derive are very similar to the ones obtained in various forms and under a variety of assumptions by Meyn and Tweedie, summarized in their book [32].

THEOREM 5. *Suppose the Markov process $(\Xi(t) : t \geq 0)$ possesses a stationary distribution $\pi$ and suppose $\Phi$ is a geometric Lyapunov function with parameters $\gamma, t_0, K$. Then,*

$$\mathbb{E}_\pi[\Phi(\Xi(0))] \leq \frac{\phi(t_0) K}{1 - \gamma}. \qquad (27)$$

Note that, in the theorem above, it is not assumed that the stationary distribution $\pi$ is unique. The claimed bounds hold for *every* stationary distribution $\pi$. The bound of course is only meaningful when $\phi(t_0)$ is finite.

We now apply Theorem 5 to establish bounds on the tail of the stationary distribution.

THEOREM 6. *Let $(\Xi(t) : t \geq 0)$ be a Markov process with a stationary distribution $\pi$, and suppose $\Phi$ is a Lyapunov function with parameters $\gamma, t_0, K$. Assume, in addition, that there exists $\theta > 0$ such that*

$$\theta L_2(\theta, t_0) \leq \gamma. \tag{28}$$

*Then $\exp(\theta \Phi(\cdot))$ is a geometric Lyapunov function with geometric drift size parameter $1 - \gamma \theta / 2$, drift time parameter $t_0$ and exception parameter $e^{\theta K}$. Consequently, for every $s > K$,*

$$\mathbb{P}_\pi(\Phi(\Xi(0)) > s) \leq (1 - \gamma \theta / 2)^{-1} L_1(\theta, t_0) \exp(-\theta(s - K)). \tag{29}$$

REMARK 6. At first glance, the division by $1 - \gamma \theta / 2$ above may seem erroneous since, in principle, this expression can be negative. We now show that condition (28) rules this out. From the definition of $L_2(\cdot, \cdot)$, we have that

$$\begin{aligned}
L_2(\theta, t_0) &\geq \sup_{x \in \mathcal{X}} \mathbb{E}_x[(\Phi(\Xi(t_0)) - \Phi(x))^2] \\
&\overset{(a)}{\geq} \sup_{x \in \mathcal{X}} (\mathbb{E}_x[\Phi(\Xi(t_0)) - \Phi(x)])^2 \\
&\overset{(b)}{\geq} \gamma^2,
\end{aligned}$$

where (a) follows from Jensen's inequality and (b) follows from the definition of $\gamma$, in particular,

$$\sup_{x : \Phi(x) > K} \mathbb{E}_x[\Phi(\Xi(t_0)) - \Phi(x)] \leq -\gamma.$$

The condition (28) then implies that $\gamma \theta \leq 1$.

REMARK 7. The upper bound in (29) can be optimized by selecting $\theta > 0$ to be the largest possible value satisfying (28). Of course, if the largest such $\theta$ is equal to zero, then the bound is trivial.

**4. Main results.** Throughout this section we consider a sequence of generalized Jackson networks $\Xi^n$ in heavy traffic, and the associated Markov process $\bar{Q}^n(t) = (Q^n(t), \hat{a}^n(t), \hat{v}^n(t))$ introduced in Section 2; the latter two coordinates denote the vectors of residual interarrival times and service times, respectively.

4.1. *Tightness of stationary distributions.* In this section we establish tightness of the sequence of stationary distributions corresponding to the normalized state process $\{Q^n(nt)/\sqrt{n}\}$. Key to this result are the tail bounds derived for general Markov processes in Section 3. To appeal to these results, we consider the workload function defined in (15), whose key properties are given in Proposition 1. The next result establishes that this workload function has the requisite attributes to serve as a Lyapunov function for the sequence of GJN's in heavy-traffic.

PROPOSITION 2. *There exist constants $t_0, c_0 > 0$ which depend only on the parameters of the network $\Xi$ and are independent of $n$, such that, for all sufficiently large, $n \geq 1$,*

$$\sup_{(z,a,v):\, w'z > c_0\sqrt{n}} \{\mathbb{E}[w'Q^n(nt_0)|\bar{Q}^n(0) = (z,a,v)] - w'z\} \leq -\sqrt{n}. \tag{30}$$

*In addition, there exists $\theta_0 > 0$ which depends only on $\Xi$ such that*

$$\limsup_{n \to \infty} \sup_{(z,a,v) \in \mathcal{X}} \mathbb{E}[\exp(n^{-1/2}\theta_0(w'Q^n(nt_0) - w'z)^+)|\bar{Q}^n(0) = (z,a,v)]$$
$$< \infty \tag{31}$$

*and*

$$\limsup_{n \to \infty} \sup_{(z,a,v) \in \mathcal{X}} n^{-1}\mathbb{E}[(w'Q^n(nt_0) - w'z)^2$$
$$\times \exp(n^{-1/2}\theta_0(w'Q^n(nt_0) - w'z)^+)|\bar{Q}^n(0) = (z,a,v)] \tag{32}$$
$$< \infty.$$

The immediate implication of the above result is the following.

PROPOSITION 3. *For all sufficiently large $n$, the function $\Phi(z,a,v) = w'z$ is a Lyapunov function with drift size parameter $-\sqrt{n}$, drift time parameter $nt_0$ and exception parameter $c_0\sqrt{n}$. In addition,*

$$\limsup_{n \to \infty} L_1(\theta_0/\sqrt{n}, nt_0) < \infty, \tag{33}$$

$$\limsup_{n \to \infty} \frac{1}{n} L_2(\theta_0/\sqrt{n}, nt_0) < \infty. \tag{34}$$

Armed with this result, we now use Theorem 6 to establish the sought tightness of the sequence of stationary distributions corresponding to the sequence of normalized GJN's.

THEOREM 7. *There exist constants $C_1, c_1 > 0$ which depend only on $\Xi$ such that, for, all sufficiently large $n$, the sequence of stationary distributions $\pi^n$ of the networks $\{\Xi^n\}$ satisfies*

$$\mathbb{P}_{\pi^n}(n^{-1/2}w'Q^n(0) > s) \leq C_1 \exp(-c_1 s), \tag{35}$$

*for all $s > 0$.*

Since $1 - \rho^{*n} = \kappa^*/\sqrt{n}$, the theorem implies that each individual queue $Q_j^n(t)$ behaves in steady state as $O(\frac{1}{1-\rho^{*n}})$ as $\rho^{*n} \to 1$.

Recall that a sequence of real valued random variables $X_n, n \in \mathbb{N}$, is said to be tight if, for every $\varepsilon > 0$, there exists sufficiently large constant $C > 0$ such that $\sup_n \mathbb{P}(|X_n| > C) < \varepsilon$. Refer to [6] for more detailed definitions and applications of this important concept. The following corollary is a consequence of the theorem above.

COROLLARY 1. *Let $\bar{Q}^n(0) = (Q^n(0), \hat{a}^n(0), \hat{v}^n(0))$ be distributed according to $\pi^n$. Then the sequence of $\mathbb{R}^J$-valued random vectors $\{Q^n(0)/\sqrt{n}\}$ is tight.*

4.2. *Interchange of limits.* We are now in position to prove the sought interchange of limits indicated in the schematic diagram in Figure 1. Let $\hat{\pi}^n$ denote the stationary distribution of the rescaled process $\hat{Q}^n(t) = Q^n(nt)/\sqrt{n}$, where, as above, $(Q^n(t): t \geq 0)$ denotes the queue length process for the GJN with traffic intensity $\rho^n$. With respect to $\pi^n$, we have, by Theorem 2, that, under some mild assumptions for each fixed $n$, the distribution of $Q^n(t)$ converges weakly to $\pi^n$ as $t \to \infty$. On the other hand, Reiman's diffusion approximation stated in Theorem 4 asserts that starting with the original GJN and taking heavy-traffic limits, that is, letting $\rho^n \uparrow 1$ as $n \to \infty$, the normalized queue length process $(\hat{Q}^n(t): t \geq 0)$ converges weakly to an $(\beta, \Gamma, I - P)$-RBM $(Z(t): t \geq 0)$ identified in Section 2. Subsequently, letting $t \to \infty$, one obtains a steady-state approximation for the original network given by the distribution $\pi_{\text{RBM}}$. We now prove that the above limits can be interchanged, that is, first taking the limit as $t \to \infty$ and letting the GJN reach steady-state, then taking heavy-traffic limits; the resulting process converges to a weak limit and this limit is $\pi_{\text{RBM}}$. This validates the steady-state heavy-traffic approximation which is built on the distribution of the approximating RBM.

To prove the interchange of limits recall the assertion of Corollary 1. In the space of probability measures on $\mathbb{R}^J$ endowed with topology of weak convergence, the tightness result of this corollary implies that there exists at least one limit point, $\pi^*$, for the sequence $\{\hat{\pi}^n\}$, and this limit is a probability measure. Our main result stated below asserts that such a limit is unique and is identical to the unique stationary distribution of the associated RBM.

THEOREM 8. *The sequence $\{\hat{\pi}^n\}$ converges weakly to $\pi_{\text{RBM}}$, where $\pi_{\text{RBM}}$ is the unique stationary distribution of the $(\beta, \Gamma, I - P)$-RBM.*

Thus, we have that the probability distribution of $\{Q^n(0)/\sqrt{n}\}$ converges weakly to $\pi_{\text{RBM}}$ when $Q^n(0)$ is distributed according to $\pi^n$. This establishes the validity of the commuting diagram depicted in Figure 1.

Given our assumptions on the primitives, in particular, the existence of a (conditional) moment generating function in the neighborhood of zero, it stands to reason that the queue length process under diffusion scaling should be uniformly integrable. Moreover, we expect the interchange argument to hold for all moments of the normalized queue length process. The following result asserts that this is indeed the case.

COROLLARY 2. *There exists a vector $\theta \in \mathbb{R}_+^J$ with all coordinates being strictly positive, such that the family of random variables $\{\exp(\theta' Q^n(0)/\sqrt{n}\,)\}$ is uniformly integrable, where $Q^n(0) \sim \pi^n$ for each $n = 1, 2, \ldots$. Consequently, for all $p \geq 1$,*

$$\mathbb{E}_{\pi^n}[n^{-1/2} Q^n(0)]^p \to \mathbb{E}_{\pi_{\mathrm{RBM}}}[Z(0)]^p \qquad as\ n \to \infty,$$

*where $\pi_{\mathrm{RBM}}$ is the stationary distribution of the RBM $(Z(t): t \geq 0)$ with parameters $(\beta, \Gamma, I - P')$.*

We note that the above result is particularly useful given that current numerical algorithms for analysis of the steady-state of RBM, see [13] and [36], are mostly effective in the computation of stationary moments (as opposed to the full distribution).

4.3. *Extension to sojourn times.* Given a job which arrives externally into station $j \in \mathcal{J}$, we will say that this job has the visit vector $h \in \mathbb{Z}_+^J$ if it visits station $i$ precisely $h_i$ times before leaving the network. Of course, not every visit vector is feasible. That is, there may not exist a job with visit vector $h$ for all $h \in \mathbb{Z}_+^J$. For example, in a tandem queueing system with two stations and no feedback, only $h = (1, 1)$ is feasible. As in [33], we let $D_{j,h}(t)$ be the sojourn time of the first job which arrives into station $j$ after time $t$ and has a visit vector $h$. Given a state $\bar{Q}(t) = (z, a, v) \in \mathcal{X}$ of the Markov process at time $t$, the distribution of $D_{j,h}(t)$ is independent from $\bar{Q}(t'), t' < t$. Thus, we may consider the random variable $D_{j,h}(z, a, v)$ which is a function of the state $(z, a, v)$. The following result is then an immediate consequence of the above observation.

PROPOSITION 4. *Let $\pi$ be the stationary distribution of the Markov process $\bar{Q}(t)$. Then for every station $j$ and feasible visit vector $h \in \mathbb{Z}_+^J$, the process $D_{j,h}(t) := D_{j,h}(\bar{Q}(t))$ is stationary.*

In the queueing system $\Xi^n$ let $D_{j,h}^n(t)$ denote $D_{j,h}^n(\bar{Q}^n(t))$. The weak convergence of the processes $D_h^n(nt)/\sqrt{n}$ defined in the queueing network $\Xi^n$ to an RBM was established by Reiman [33]. The result is stated there under the assumption that the queueing network is initially empty $Q^n(0) = 0$. The techniques used by Reiman allow one to obtain a stronger result, which we will use in the present paper. We do not give a full proof of this stronger version, but just highlight the necessary adjustments to the basic arguments given in [33].

THEOREM 9. *Consider a sequence of GJN, $\{\Xi^n\}$, in heavy traffic. Suppose $Q^n(0)/\sqrt{n}$ converges weakly to some limiting distribution $\pi_0$. Then for every station $j$ and every feasible $h$, the process $(D^n_{j,h}(t')/\sqrt{n} : 0 \le t' \le nt)$ converges weakly to $(h'MZ(t') : 0 \le t' < \infty)$, where $Z(t)$ is the $(\beta, \Gamma, I - P')$-RBM with $Z \sim \pi_0$.*

We now state our main result for the sojourn times in steady-state. The proof is immediate by combining Theorems 8 and 9.

THEOREM 10. *Consider the sequence of GJN, $\{\Xi^n\}$, in heavy-traffic. Suppose $\bar{Q}^n(0) \sim \pi^n$. Then for every station $j$ and every feasible $h \in \mathbb{Z}^J_+$, the convergence*

$$\frac{D^n_{j,h}(0)}{\sqrt{n}} \Rightarrow h'(\mu_1^{-1} Z_1(0), \ldots, \mu_J^{-1} Z_J(0))$$

*holds, where $Z(0) \sim \pi_{\text{RBM}}$.*

4.4. *Product form RBM.* As established in [19], in special cases described below the stationary distribution of RBM has an explicit product form exponential distribution. Our main result, Theorem 8, can then be used to validate closed form steady-state approximations for certain classes of GJN in heavy traffic. In this subsection we illustrate this application of our main theorem.

The following result follows from a more general theorem on RBM in steady state, established by Harrison and Williams [19]; see also [11], Theorem 7.12.

THEOREM 11 (Product form RBM stationary distribution). *Given a $(\beta, \Gamma, I - P')$-RBM, suppose $[I - P']^{-1}\beta < 0$ and, in addition,*

(36) $$2\Gamma = [I - P']D^{-1}\Lambda + \Lambda D^{-1}[I - P],$$

*where $D$ is the diagonal matrix with elements $(1 - p_{jj}), 1 \le j \le J$, and $\Lambda$ is the diagonal matrix with elements $\Gamma_{jj}, 1 \le j \le J$. Then the stationary distribution $\pi_{\text{RBM}}$ has a product form density*

$$f(z_1, \ldots, z_J) = \prod_{1 \le j \le J} \eta_j \exp(-\eta_j z_j), \qquad (z_1, \ldots, z_J) \in \mathbb{R}^J_+,$$

*where $\eta = -2\Lambda^{-1} D[I - P']^{-1}\beta$.*

It can be checked that condition (36), which is also referred to as the *skew symmetry condition*, holds when $c_{a,j}^2 = c_{s,j}^2 = 1$. In particular, this includes the classical Jackson network case. Consider the GJN, $\{\Xi^n\}$, in heavy traffic, and recall that $\beta = -[I - P']M^{-1}\kappa$ with $\kappa = (\kappa_1, \ldots, \kappa_J)$. Using this, we have $\eta_j = 2\Gamma_{jj}^{-1}(1 - p_{jj})\mu_j\kappa_j$. Combining with Theorem 8, we obtain the following theorem, which gives an explicit expression for the limiting steady-state distribution of GJN in heavy traffic.

THEOREM 12. *Consider a sequence of GJN, $\{\Xi^n\}$, in heavy traffic. Suppose condition (36) holds. Then the sequence of stationary measures $\pi^n$ satisfies*

$$\lim_{n \to \infty} \mathbb{P}_{\pi^n}\left(\frac{Q^n(0)}{\sqrt{n}} > z\right) = \prod_{1 \leq j \leq J} \exp(-2\Gamma_{jj}^{-1}(1 - p_{jj})\mu_j\kappa_j z_j),$$

*for every $z \in \mathbb{R}_+^J$.*

Recalling that $\rho_j = 1 - \kappa_j/\sqrt{n}$, we can informally express the assertion of the above limit theorem as

$$\mathbb{P}_\pi((1 - \rho_j)Q_j(0) > z_j, \forall j) \approx \prod_{1 \leq j \leq J} \exp(-2\Gamma_{jj}^{-1}(1 - p_{jj})\mu_j z_j),$$

when $\rho_j \approx 1$, for all $j = 1, \ldots, J$, where $\pi$ is a stationary distribution of GJN. Thus, for certain classes of GJN for which the skew symmetry condition holds, we obtain a simple approximation to the steady-state on the basis of the product form RBM stationary distribution.

**5. Concluding remarks.** We considered a generalized Jackson network (GJN) in heavy traffic and proved that the stationary distribution of the GJN under diffusion scaling converges to the stationary distribution of the associated reflected Brownian motion, thus resolving the open problem of so called "interchange of limits" for this class of open queueing networks. There are several interesting questions which remain unresolved and are worth pursuing in future research.

1. The bounds on the stationary distribution of the queue lengths in GJN obtained in the present paper are not explicit, that is, their dependence on the stochastic primitives is not articulated in full detail. The techniques used in the paper can be applied, in principle, to obtain more explicit expressions. We are currently pursuing such results in ongoing work, deriving explicit upper and lower bounds on the steady-state distribution of a GJN.

2. The assumptions adopted in the present paper are certainly not the tightest possible. In particular, it should be possible to carry out a similar analysis when one only assumes finite $2 + \delta$ moments, for any $\delta > 0$, for the interarrival and service times. (These conditions are very close to being necessary and sufficient for the functional central limit theorem to hold, and thus for proving heavy-traffic process limits.) The Lyapunov approach we develop in the current paper can be adjusted to provide polynomial tail bounds under such assumptions; in essence, a $p$th moment of the primitives would yield existence of the $(p-1)$st moment of the stationary queue length (see, e.g., [14]). This would be used in conjunction with strong approximation coupling inequalities that are available under polynomial moment conditions on the primitives, the weakest being the existence of a $2+\delta$ moment. In the present paper we assume exponential moments and consequently obtain exponential tail bounds on the steady-state distribution of the GJN.

3. The main results in our current paper hinge on the use of a simple linear Lyapunov function. There is a long history of applying Lyapunov function methods for purposes of performance analysis, specifically for multiclass queueing networks; see [2, 3, 4, 5, 26, 29, 30]. It should be noted, however, that in almost all the existing literature, Markovian-type queueing networks are considered. The extension of these bounds to general type queueing networks has not been explored to date.

4. This paper presents what is, in our view, a very interesting link between diffusion approximations of queueing networks and Lyapunov function methods. The latter have been traditionally used to establish stability and obtain performance bounds for queueing networks, as well as diffusion limits that correspond to the underlying queueing model (see [16] for an example of the latter). To the best of our knowledge, such methods have not yet been applied to queueing networks to establish properties of their diffusion approximations (e.g., to ascertain tightness of the sequence of steady-state distributions, as is established in the current paper). The proof technique introduced in this paper opens up the possibility of establishing similar results for other models where: (i) process-level diffusion approximations have been established; (ii) the stationary distribution of both the underlying queueing model and the associated diffusion limit is known to exist; and no rigorous link between (i) and (ii) is known. Notable examples of such systems are feedforward (acyclic) queueing networks and certain instances of multiclass networks for which stability conditions are known and Brownian approximations have been rigorously established. One difficulty in directly applying the methods of this paper is the apparent lack of strong solution and Lipschitz continuity in Skorohod mappings corresponding to these models. Proving an interchange of limits for these classes of queueing networks is an interesting direction of future research that may benefit from the results we develop in the current paper.

## APPENDIX: PROOFS

### A.1. Proofs of the main results.

PROOF OF THEOREM 5. By definition of a geometric Lyapunov function (23), we have that, for every $x \in \mathcal{X}$,

$$\mathbb{E}_x[\Phi(\Xi(t_0))] \leq \max\{\gamma\Phi(x), \phi(t_0)\Phi(x)\mathbb{1}\{\Phi(x) \leq K\}\}$$
$$\leq \gamma\Phi(x) + \phi(t_0)K.$$

Fix $k \in \mathbb{N}$ and put $\Phi_k(x) := \Phi(x) \wedge k$. Then, if $\Phi_k(x) \leq k$

$$\Phi_k(x) - \mathbb{E}_x[\Phi_k(\Xi(t_0))] = \Phi(x) - \mathbb{E}_x[\Phi_k(\Xi(t_0))]$$
$$\geq \Phi(x) - \mathbb{E}_x[\Phi(\Xi(t_0))]$$
$$\geq (1-\gamma)\Phi(x) - \phi(t_0)K$$
$$\geq -\phi(t_0)K.$$

On the other hand, if $\Phi_k(x) > k$, then

$$\Phi_k(x) - \mathbb{E}_x[\Phi_k(\Xi(t_0))] \geq k - \mathbb{E}_x[\Phi_k(\Xi(t_0))] \geq 0.$$

Thus, $\Phi_k(x) - \mathbb{E}_x[\Phi_k(\Xi(t_0))]$ is bounded from below by $-\phi(t_0)K$ for all $x \in \mathcal{X}$. It is also bounded from above by $k$. Now, the monotone convergence theorem implies that $\mathbb{E}_x\Phi_k(\Xi(t_0)) \uparrow \mathbb{E}_x\Phi(\Xi(t_0))$ as $k \to \infty$. So, $\Phi_k(x) - \mathbb{E}_x[\Phi_k(\Xi(t_0))] \to \Phi(x) - \mathbb{E}_x[\Phi(\Xi(t_0))]$ as $k \to \infty$, for all $x \in \mathcal{X}$. Since $\Phi_k(x) - \mathbb{E}_x[\Phi_k(\Xi(t_0))]$ bounded below, we can appeal to Fatou's lemma to conclude that

(37)
$$\int (\Phi(x) - \mathbb{E}_x[\Phi(\Xi(t_0))])\pi(dx)$$
$$\leq \liminf_{k \to \infty} \int (\Phi_k(x) - \mathbb{E}_x[\Phi_k(\Xi(t_0))])\pi(dx) = 0,$$

where the equality follows from stationarity of $\pi$. On the other hand, by definition of the Lyapunov function, we have that

$$\Phi(x) - \mathbb{E}_x[\Phi(\Xi(t_0))] \geq (1-\gamma)\Phi(x) - \phi(t_0)K.$$

Integrating both sides with respect to $\pi$ and using (37), we conclude that

$$\mathbb{E}_\pi \Phi(\Xi(0)) \leq \frac{\phi(t_0)K}{1-\gamma},$$

which concludes the proof. □

PROOF OF THEOREM 6. Let $\theta > 0$ satisfy (28). Introduce the function $\Upsilon : \mathcal{X}^2 \times \mathbb{R}_+ \to \mathbb{R}_+$, defined as

$$\Upsilon(x_1, x_2, y) := \exp(y(\Phi(x_2) - \Phi(x_1))).$$

Using a second-order Taylor's expansion around $y = 0$, we have, for some $0 < \theta' \leq \theta$,

$$\Upsilon(x_1, x_2, \theta) = 1 + \theta(\Phi(x_2) - \Phi(x_1))$$
$$+ \frac{\theta^2}{2}(\Phi(x_2) - \Phi(x_1))^2 \exp(\theta'(\Phi(x_2) - \Phi(x_1))),$$
$$\leq 1 + \theta(\Phi(x_2) - \Phi(x_1))$$
$$+ \frac{\theta^2}{2}(\Phi(x_2) - \Phi(x_1))^2 \exp(\theta(\Phi(x_2) - \Phi(x_1))^+).$$

We now fix arbitrary $x$ such that $\Phi(x) > K$, set $x_1 = x, x_2 = \Xi(t_0)$ and take the expectation of both sides above to obtain

$$\mathbb{E}_x[\Upsilon(x, \Xi(t_0), \theta)]$$
$$\leq 1 + \theta \mathbb{E}_x[\Phi(\Xi(t_0)) - \Phi(x)]$$
$$+ \frac{\theta^2}{2} \mathbb{E}_x[(\Phi(\Xi(t_0)) - \Phi(x))^2 \exp(\theta(\Phi(\Xi(t_0)) - \Phi(x))^+)]$$
$$\stackrel{(a)}{\leq} 1 - \gamma\theta + \frac{\theta^2}{2} L_2(\theta, t_0)$$
$$\stackrel{(b)}{\leq} 1 - \gamma\theta/2,$$

where (a) follows from (22) and the definition of $L_2$, and (b) follows from (28). Note that the condition $\Phi(x) > K$ is equivalent to $\exp(\theta\Phi(x)) > \exp(\theta K)$. Thus, $\exp(\theta\Phi(x))$ is a geometric Lyapunov function with parameters $\gamma\theta/2$, $t_0, \exp(\theta K)$. This concludes the proof of the first part of Theorem 6.

To prove the second part, observe that, for the constructed Lyapunov function,

(38) $\qquad \phi(t_0) = \sup_x \mathbb{E}_x[\exp(\theta(\Phi(\Xi(t_0)) - \Phi(x)))] = L_1(\theta, t_0).$

We use Markov's inequality, (38) and the bound (27) of Theorem 5, which yields

$$\mathbb{P}_\pi(\Phi(\Xi(0)) > s) = \mathbb{P}_\pi(\exp(\theta\Phi(\Xi(0))) > \exp(\theta s))$$
$$\leq \exp(-\theta s)\mathbb{E}_\pi[\exp(\theta\Phi(\Xi(0)))]$$
$$\leq \exp(-\theta s)\frac{L_1(\theta, t_0)\exp(\theta K)}{1 - \gamma\theta/2}$$
$$= (1 - \gamma\theta/2)^{-1} L_1(\theta, t_0) \exp(-\theta(s - K)).$$

This completes the proof of the second part of the theorem. □

PROOF OF PROPOSITION 2.  The key to the proof is the following result.

LEMMA A.1. *Let $(X^n(t) : t \geq 0)$ be the net input process for the GJN, $\Xi^n$, and let $(x_z^n(t) : t \geq 0)$ denote its fluid counterpart starting at state $z$ [i.e., $x_z^n(0) = z$]. Then, there exists $\theta_1 > 0$ such that*

$$\text{(39)} \quad \limsup_{n \to \infty} \sup_{(z,a,v) \in \mathcal{X}} n^{-1/2} \mathbb{E}\left[\sup_{0 \leq t \leq n} \|X^n(t) - x_z^n(t)\| \,\Big|\, \bar{Q}^n(0) = (z,a,v)\right] < \infty,$$

$$\text{(40)} \quad \limsup_{n \to \infty} \sup_{(z,a,v) \in \mathcal{X}} \mathbb{E}\left[\sup_{0 \leq t \leq n} \exp(n^{-1/2} \theta_1 \|X^n(t) - x_z^n(t)\|) \,\Big|\, \bar{Q}^n(0) = (z,a,v)\right] < \infty$$

*and*

$$\text{(41)} \quad \limsup_{n \to \infty} \sup_{(z,a,v) \in \mathcal{X}} n^{-1} \mathbb{E}\bigg[\sup_{0 \leq t \leq n} \|X^n(t) - x_z^n(t)\|^2 \\ \times \exp(n^{-1/2} \theta_1 \|X^n(t) - x_z^n(t)\|) \,\Big|\, \bar{Q}^n(0) = (z,a,v)\bigg] < \infty.$$

The proof of this lemma is given in Section A.2.

From (8), (9) and (10), it follows that $(Y^n(\cdot), Q^n(\cdot))$ jointly solve the Skorohod problem for the unreflected net input process $X^n(\cdot)$. That is, $Y^n(\cdot) = \Psi_1(X^n(\cdot)), Q^n(\cdot) = \Psi_2(X^n(\cdot))$. Let $(y^n, q^n) = (\Psi_1(x_z^n), \Psi_2(x_z^n))$ be the solution of the Skorohod problem for the linear function $x_z^n(\cdot)$. Fix $\theta > 0$ and $t_0 > 0$. Then,

$$n^{-1/2} \theta (w' Q^n(t_0 n) - w' z)^+$$
$$\leq n^{-1/2} \theta (w' Q^n(t_0 n) - w' q^n(t_0 n))^+ + n^{-1/2} \theta (w' q^n(t_0 n) - w' z)^+$$
$$\text{(42)} \quad \overset{(a)}{\leq} n^{-1/2} \theta \sup_{0 \leq t \leq t_0 n} |w' Q^n(t) - w' q^n(t)|$$
$$\overset{(b)}{\leq} C_1 \theta n^{-1/2} \sup_{0 \leq t \leq t_0 n} \|Q^n(t) - q^n(t)\|$$
$$\overset{(c)}{\leq} R_1 C_1 \theta n^{-1/2} \sup_{0 \leq t \leq t_0 n} \|X^n(t) - x_z^n(t)\|.$$

In the above (a) follows from from Proposition 1, from which we have that, for every $z \in \mathbb{R}_+^J$,

$$\text{(43)} \quad w' q^n(t_0 n) \leq \left(w' z - \min_{1 \leq j \leq J} \mu_j (1 - \rho_j^n) t_0 n\right)^+$$

$$= \left(w'z - \min_{1 \le j \le J} \mu_j t_0 \kappa^* \sqrt{n}\right)^+.$$

Thus, the second term on the RHS of the first inequality above is identically zero. The inequality (b) follows from setting $C_1 = \max w_j$, and (c) follows from (13) of Theorem 1. Thus, setting $\theta_0 = \theta_1/(R_1 C_1 \sqrt{t_0})$, where $\theta_1$ is identified in Lemma A.1, concludes the proof of (31) and (32) for arbitrary choice of $t_0$.

We now concentrate on proving (30). From (42) and (39) of Lemma A.1, we have that there exists a constant $C_2 > 0$ such that

$$\mathbb{E}[|w'Q^n(t_0 n) - w'q^n(t_0 n)||\bar{Q}^n(0) = (z, a, v)]$$
$$\le \mathbb{E}\left[\sup_{0 \le t \le t_0 n} |w'Q^n(t) - w'q^n(t)||\bar{Q}^n(0) = (z, a, v)\right]$$
$$\le C_2 \sqrt{t_0 n},$$

holds for all $(z, a, v) \in \mathcal{X}$, $t_0 > 0$ and all sufficiently large $n$, where we have substituted $t_0 n$ for $n$ in (39). Combining the above with (43), we obtain that, for every $t_0 > 0$, $(z, a, v) \in \mathcal{X}$, and all sufficiently large $n$,

$$\mathbb{E}[w'Q^n(t_0 n)|\bar{Q}^n(0) = (z, a, v)]$$
$$\le \max\left\{C_2 \sqrt{t_0 n}, C_2 \sqrt{t_0 n} + w'z - \min_{1 \le j \le J} \mu_j t_0 \kappa^* \sqrt{n}\right\}$$
$$= w'z + \max\left\{C_2 \sqrt{t_0 n} - w'z, C_2 \sqrt{t_0 n} - \min_{1 \le j \le J} \mu_j t_0 \kappa^* \sqrt{n}\right\}.$$

We now set $t_0$ as a function of $C_2, \kappa^*$ and $\min_j\{\mu_j\}$, so that $C_2 \sqrt{t_0} - \min_{1 \le j \le J} \mu_j t_0 \kappa^* \le -1$. This is always feasible for large enough value of $t_0$ since the expression above is quadratic in $\sqrt{t_0}$ with negative coefficient in front of $(\sqrt{t_0})^2$. For said $t_0$, we obtain $C_2 \sqrt{t_0 n} - \min_{1 \le j \le J} \mu_j t_0 \kappa^* \sqrt{n} \le -\sqrt{n}$. Then we set $c_0 = C_2 \sqrt{t_0} + 1$. Consequently, for values of $z$ such that $w'z > c_0 \sqrt{n}$, we have that $C_2 \sqrt{t_0 n} - w'z \le -\sqrt{n}$. We conclude that for our constants $t_0, c_0$ and $C_2$, we have that

$$\mathbb{E}[w'Q^n(t_0 n)|\bar{Q}^n(0) = (z, a, v)] \le w'z - \sqrt{n}$$

whenever $w'z > c_0 \sqrt{n}$. This proves (30). $\square$

PROOF OF THEOREM 7. We apply Proposition 3. Fix a constant $0 < c < 1$. The proposition implies that, for $\theta_n \equiv c\theta_0 n^{-1/2}$, we have

(44) $$\limsup_{n \to \infty} L_1(\theta_n, nt_0) < \infty$$

and

(45) $$\limsup_{n \to \infty} n^{-1} L_2(\theta_n, nt_0) < \infty.$$

This implies that
$$\limsup_{n\to\infty} n^{-1/2}\theta_n L_2(\theta_n, nt_0) = \limsup_{n\to\infty} c\theta_0 n^{-1} L_2(\theta_n, nt_0) \leq 1,$$
provided that $c$ is sufficiently small. This means that condition (28) of Theorem 6 holds for all sufficiently large $n$.

Applying Theorem 6 and using $n^{1/2}\theta_n = c\theta_0$, we obtain, for every $s > 0$,
$$\mathbb{P}_{\pi_n}(w'Q^n(0)/\sqrt{n} > s) \leq \frac{L_1(\theta_n, nt_0)}{1 - c\theta_0/2} \exp(-\theta_n(sn^{1/2} - c_0 n^{1/2}))$$
$$= \frac{L_1(\theta_n, nt_0)}{1 - c\theta_0/2} \exp(-c\theta_0(s - c_0)).$$

From (44), there exists a sufficiently large $n_0$ and a constant $c_3 > 0$ such that $L_1(\theta_n, nt_0) < c_3$ for all $n > n_0$. Then for all $s > c_0$ and all $n > n_0$,
$$\mathbb{P}_{\pi_n}(w'Q^n(0)/\sqrt{n} > s) \leq \frac{c_3}{1 - c\theta_0/2} \exp(-c\theta_0(s - c_0)).$$

Setting $C_1$ and $c_1$ appropriately, we obtain the result. □

PROOF OF THEOREM 8. Corollary 1 implies that, for every subsequence $\hat{\pi}^{n_k}$, there exists a weakly convergent subsequence $\hat{\pi}^{n_{k_i}}$. Let $\pi^*$ denote any weak limit of this subsequence. We claim that every such limit measure $\pi^*$ is equal to $\pi_{\text{RBM}}$, implying that all such weak limits are the same. It then follows that $\hat{\pi}^n \Rightarrow \pi_{\text{RBM}}$ as $n \to \infty$. Thus, we are left to show that every converging sequence $\hat{\pi}^{n_{k_i}} \Rightarrow \pi^*$ converges to $\pi_{\text{RBM}}$. For notational simplicity, we denote $\hat{\pi}^{n_{k_i}}$ simply by $\hat{\pi}^n$.

We apply Theorem 4. Fix an arbitrary $t > 0$. Since $\hat{\pi}^n$ converges weakly to $\pi^*$, then the process $(Q^n(nt')/\sqrt{n}: 0 \leq t' \leq t)$ conditioned on $Q^n(0) \sim \pi^n$ converges weakly to a $(\beta, \Gamma, I - P')$-RBM $(Z(t'): 0 \leq t' \leq t)$ with initial distribution given by $\pi^*$. In particular, $Q^n(nt)/\sqrt{n}$ converges weakly to $Z(t)$ with $Z(0)$ distributed according to $\pi^*$. But from stationarity of $\pi^n$, $Q^n(nt)/\sqrt{n}$ is distributed as $\hat{\pi}^n$, which, by assumption, converges to $\pi^*$. Therefore, $Z(t)$ is distributed as $\pi^*$ and $\pi^*$ must be the unique stationary distribution: $\pi^* = \pi_{\text{RBM}}$. □

PROOF OF COROLLARY 2. It suffices to prove the asserted uniform integrability. Fix $\delta > 0$. By Theorem 7, we have that, for all sufficiently large $n$,
$$\mathbb{E}[\exp(n^{-1/2}\delta w'Q^n(0))] = \int_0^\infty \mathbb{P}(n^{-1/2}w'Q^n(0) > \delta^{-1}\log x)\,dx$$
$$\leq \int_0^\infty \exp(-c_1\delta^{-1}\log x + c_2)\,dx < \infty,$$
provided that $c_1\delta^{-1} > 1$. Setting $\theta = \delta w$, we obtain the result. □

PROOF OF THEOREM 9. The proof repeats the argument used in [33]. It is first shown in [33] that, when $Q^n(0) = 0$,

$$\lim_{n \to \infty} \mathbb{E}\bigg[\sup_{0 \le t' \le nt} D_{j,h}^n(t')/n\bigg] = 0, \tag{46}$$

for every $j \in \mathcal{J}$. This convergence is shown by bounding the supremum by a certain linear functional of $\sup_{0 \le t' \le nt} \|Q^n(t')\|/n$. Then a random time change theorem is applied to show that the distribution of sojourn times in station $j$ at time $nt'$ becomes "indistinguishable" from the distribution of the sum of the workloads in every individual station $j'$ multiplied by $h_i$. The latter converges weakly to $h_{j'} \mu_{j'} Z_j(t')$ by Theorem 4, and the result follows. But (46) also holds in our setting when $Q^n(0)$ is not zero, but instead $Q^n(0)/\sqrt{n} \Rightarrow \pi_0^*$. Indeed, by Theorem 4, $(Q^n(t')/\sqrt{n} : 0 \le t' \le nt)$ converges weakly to a corresponding RBM with $Z(0) \sim \pi_0^*$. Then $\sup_{0 \le t' \le nt} \|Q^n(t')\|/\sqrt{n}$ converges weakly to $\sup_{0 \le t' \le t} \|Z(t')\|$, with initial distribution $\pi_0^*$. This implies that, for every $\varepsilon > 0$,

$$\lim_{n \to \infty} \mathbb{P}\bigg(\sup_{0 \le t' \le nt} \|Q^n(t')\|/n > \varepsilon\bigg) = \lim_{n \to \infty} \mathbb{P}\bigg(\sup_{0 \le t' \le nt} \|Q^n(t')\|/\sqrt{n} > \varepsilon \sqrt{n}\bigg) = 0,$$

since $\mathbb{E}_{\pi_0^*}[\sup_{0 \le t' \le t} \|Z(t')\|] < \infty$. Repeating the argument in [33], (46) then holds as well and the statement in the theorem follows. $\square$

### A.2. Proofs of auxiliary results.

PROOF OF LEMMA A.1. We first show that (40) implies (39) and (41). Given $\theta_1 > 0$, select $C$ large enough so that $\exp(\theta_1 x) > x^2 > x$ for all $x > C$. Replacing $n^{-1/2} \sup_{0 \le t \le n} |X^n(t) - x_z^n(t)|$ with $\max(C, n^{-1/2} \sup_{0 \le t \le n} \|X^n(t) - x_z^n(t)\|)$, we obtain that (39) and (41) hold, provided that

$$\limsup_{n \to \infty} \sup_{(z,a,v) \in \mathcal{X}} \mathbb{E}\bigg[\sup_{0 \le t \le n} \exp(2\theta_1 \max(C, n^{-1/2}\|X^n(t) - x_z^n(t)\|))$$
$$|\bar{Q}^n(0) = (z,a,v)\bigg] < \infty.$$

But this bound holds if and only if (40) holds. Thus, it suffices to prove (40).

The proof of (40) relies on considering tail bounds on the uniform deviations of $X^n$ from $x_z^n$. We expand $x_z^n(t)$ as $(x_{z,1}^n(t), \ldots, x_{z,J}^n(t))$. Note that, because $X^n(0) = x_z^n(0) = Q^n(0) = z$, we have, for all $j = 1, \ldots, J$,

$$X_j^n(t) - x_{z,j}^n(t) = (A_j^n(t) - \alpha_j^n t)$$
$$+ \sum_{1 \le i \le J} p_{ij}(S_i(B_i^n(t)) - \mu_i B_i^n(t)) - (S_j(B_j^n(t)) - \mu_j B_j^n(t))$$
$$+ \sum_{1 \le i \le J} (R_j^i(S_i(B_i^n(t))) - p_{ij} S_i(B_i^n(t))),$$

where we recall that the process $S(t)$ is independent from $n$.

Fix $j, 1 \leq j \leq J$. We then have, for all integer $n \geq 1$ and $(z, a, v) \in \mathcal{X}$,

$$\mathbb{P}\left(\sup_{0 \leq t \leq n} \exp(\theta n^{-1/2}|X_j^n(t) - x_{z,j}^n(t)|) \geq u | \bar{Q}^n(0) = (z, a, v)\right)$$

$$= \mathbb{P}\left(\sup_{0 \leq t \leq n} |X_j^n(t) - x_{z,j}^n(t)| \geq \theta^{-1} n^{1/2} \log u | \bar{Q}^n(0) = (z, a, v)\right)$$

$$\leq \mathbb{P}\left(\sup_{0 \leq t \leq n} |A_j^n(t) - \alpha_j^n t| \geq (\theta^{-1} n^{1/2} \log u)/3 | \bar{Q}^n(0) = (z, a, v)\right)$$

(47)
$$+ \mathbb{P}\left(\sup_{0 \leq t \leq n} |S_j(B_j^n(t)) - \mu_j B_j^n(t)| \geq (\theta^{-1} n^{1/2} \log u)/(3(J+1))\right.$$

$$\left. | \bar{Q}^n(0) = (z, a, v)\right)$$

$$+ \mathbb{P}\left(\sup_{0 \leq t \leq n} |R_j^i(S_i(B_i^n(t))) - p_{ij} S_i(B_i^n(t))| \geq (\theta^{-1} n^{1/2} \log u)/(3J)\right.$$

$$\left. | \bar{Q}^n(0) = (z, a, v)\right),$$

for all $u \geq 1$. Since $B_j^n(t) - B_j^n(s) \leq (t - s)$ for all $t \geq s$, it suffices to bound the uniform deviation of $S_j(t)$ from $\mu_j t$, and $R_j^i(S_j(t))$ from $p_{ij} S_i(t)$ in the latter two terms on the right-hand side above, for all $i, j = 1, \ldots, J$. We now proceed to bound these terms using the so-called coupling probability bounds that are the main building blocks in strong approximations (cf. [8], Chapter 2.6). Introduce the following $J(J+2)$ constants:

(48)
$$\Gamma_j^0 = \alpha_j c_{a,j}^2 \qquad \text{for } j = 1, \ldots, J,$$
$$\Gamma_j^{J+1} = \mu_j c_{s,j}^2 \qquad \text{for } j = 1, \ldots, J,$$
$$\Gamma_{kj} = p_{kj}(1 - p_{kj}) \qquad \text{for } k, j = 1, \ldots, J.$$

We now use the following lemma, which provides bounds on the probability of deviation between the stochastic primitives and certain Brownian motion processes that are constructed on the same underlying probability space.

LEMMA A.2. *There exist $J+2$ mutually independent standard $J$-dimensional Brownian motions $W^k = (W^k(t) : t \geq 0)$, $k = 0, \ldots, J+1$, such that*

$$\mathbb{P}\left(\sup_{0 \leq t \leq n} |A_j(t) - \alpha_j t - (\Gamma_j^0)^{1/2} W_j^0(t)| \geq C_{0,j}^1 \log n + u\right.$$

$$\left. | \bar{Q}(0) = (z, 0, 0)\right) \leq C_{0,j}^2 e^{-c_{0,j} u},$$

$$\mathbb{P}\bigg(\sup_{0\leq t\leq n} |S_j(t) - \mu_j t$$
$$- (\Gamma_j^{J+1})^{1/2} W_j^{J+1}(t)| \geq C_{J+1,j}^1 \log n + u$$
$$|\bar{Q}(0) = (z,0,0)\bigg) \leq C_{J+1,j}^2 e^{-c_{J+1,j} u},$$

$$\mathbb{P}\bigg(\sup_{0\leq t\leq n} |R_j^k(\lfloor t \rfloor) - p_{kj} t - (\Gamma_{kj})^{1/2} W_j^k(t)| \geq C_{k,j}^1 \log n + u$$
$$|\bar{Q}(0) = (z,0,0)\bigg) \leq C_{k,j}^2 e^{-c_{k,j} u},$$

for $j = 1, \ldots, J$ and $k = 1, \ldots, J$, where $C_{k,j}^1$, $C_{k,j}^2$ and $c_{k,j}$, $k = 0, \ldots, J+1$, $j = 1, \ldots, J$, are finite positive constants independent of $u$ and $n$.

Note that the first identity is given for the limiting process $A(t)$ and not the indexed process $A^n(t)$. But since $A^n(t) = A(\alpha_j^{-1} \alpha_j^n t)$, $\alpha_j^n t = (\alpha_j^{-1} \alpha_j^n t) \alpha_j t$ and $\alpha_j^n < \alpha_j$, then $\sup_{0\leq t\leq n} |A_j^n(t) - \alpha_j^n t| \leq \sup_{0\leq t\leq n} |A_j(t) - \alpha_j t|$. Therefore, it suffices to obtain an appropriate bound on $\sup_{0\leq t\leq n} |A_j(t) - \alpha_j t|$. Note also that in the first expression above conditioning on $\bar{Q}(0) = (z,0,0)$ is equivalent to conditioning on $a_0 = 0$, since the arrival process is independent of the service time process and the queue length vector at time $t = 0$. Similarly, for the second expression, the conditioning is equivalent to conditioning on $v_0 = 0$; for the third expression, conditioning on $\bar{Q}(0) = (z,0,0)$ is redundant.

Using Lemma A.2, we can bound the three terms on the right-hand side of (47). Since the proof is essentially identical (with obvious modifications) in all three cases, we will provide a complete argument only for the first term on the right-hand side of (47), that is, $\mathbb{P}(\sup_{0\leq t\leq n} |A_j(t) - \alpha_j t| \geq (\theta^{-1} n^{1/2} \log u)/3 | \bar{Q}(0) = (z,a,v))$, omitting the other cases. Let $\tau = (\theta^{-1} n^{1/2} \log u)/3$. Then, we have

$$\mathbb{P}\bigg(\sup_{0\leq t\leq n} |A_j(t) - \alpha_j t| \geq \tau | \bar{Q}(0) = (z,a,v)\bigg)$$
$$= \mathbb{P}\bigg(\sup_{0\leq t\leq n} |A_j(t) - \alpha_j t| \geq \tau | \hat{a}_j(0) = a_j\bigg)$$
$$\leq \mathbb{P}\bigg(\sup_{0\leq t\leq n} |A_j(t) - \alpha_j t| \geq \tau | a_j(0) \leq \tau/(2\alpha_j), \hat{a}_j(0) = a_j\bigg)$$
(49)
$$\quad \times \mathbb{P}(a_j(0) \leq \tau/(2\alpha_j) | \hat{a}_j(0) = a_j)$$
$$\quad + \mathbb{P}(a_j(0) > \tau/(2\alpha_j) | \hat{a}_j(0) = a_j)$$

STEADY-STATE APPROXIMATIONS IN OPEN QUEUEING NETWORKS 31$$\leq \mathbb{P}\left(\sup_{0\leq t\leq n} |A_j(t) - \alpha_j t| \geq \tau | a_j(0) \leq \tau/(2\alpha_j), \hat{a}_j(0) = a_j\right)$$

$$+ \mathbb{P}(a_j(0) > \tau/(2\alpha_j)|\hat{a}_j(0) = a_j).$$

We now bound each of the two terms appearing in the sum above, starting with the second term. Using Markov's inequality, we have that

$$\mathbb{P}(a_j(0) > \tau/(2\alpha_j)|\hat{a}_j(0) = a_j) \leq \mathbb{E}[\exp(\theta a_j(0))|\hat{a}_j(0) = a_j]\exp(-\theta\tau/(2\alpha_j))$$

$$= \mathbb{E}[\exp(\theta a_j(0))|\hat{a}_j(0) = a_j]u^{-(6\alpha_j)^{-1}n^{1/2}},$$

using the definition of $\tau$. Recall that the distribution of $a_j(0)$ conditioned on $\hat{a}_j(0) = a_j$ is the same as the distribution of $a_j(1) - a_j$ conditioned on $a_j(1) > a_j$, and from (1), we have that, for $\theta < \theta^*$, $\sup_{a_j \in \mathbb{R}_+} \mathbb{E}[\exp(\theta a_j(0))|\hat{a}_j(0) = a_j] < \infty$. When $n > 2^2(6\alpha_j)^2$, we have $u^{-(6\alpha_j)^{-1}n^{1/2}} < u^{-2}$ (the choice of 2 as a lower end for the interval of integration below will become clear in what follows). This implies that, when $n > 2^2(6\alpha_j)^2$ and $\theta < \theta^*$, we have

$$\limsup_{n\to\infty} \sup_{(z,a,v)} \int_2^\infty \mathbb{P}(a_j(0) > \tau/(2\alpha_j)|\hat{a}_j(0) = a_j)\,du$$

$$(50) \quad \leq \limsup_{n\to\infty} \sup_{a_j \in \mathbb{R}_+} \mathbb{E}[\exp(\theta a_j(0))|\hat{a}_j(0) = a_j]\int_2^\infty u^{-(6\alpha_j)^{-1}n^{1/2}}\,du < \infty.$$

Now we estimate the first term in (49). For every $\tau' \leq \tau/(2\alpha_j)$, conditioning on $a_j(0) = \tau'$, we have $A_j(0,t) = 1 + A_j(\tau',t)$. Observe that the distribution of $A_j(\tau',t)$ is the same as $A_j(0, t-\tau')$ conditioned on having the residual interarrival time $a_j(0) = 0$. This follows since, by assumption, there was an arrival at time 0 (time $\tau'$ before the shift). Thus, we obtain

$$\mathbb{P}\left(\sup_{0\leq t\leq n} |A_j(t) - \alpha_j t| \geq \tau | a_j(0) = \tau', \hat{a}_j(0) = a_j\right)$$

$$= \mathbb{P}\left(\sup_{\tau'\leq t\leq n} |1 + A_j(t-\tau') - \alpha_j t| \geq \tau | \hat{a}_j(0) = 0\right)$$

$$\leq \mathbb{P}\left(\sup_{\tau'\leq t\leq n} |A_j(t-\tau') - \alpha_j(t-\tau')| + 1 + \alpha_j\tau' \geq \tau | \hat{a}_j(0) = 0\right)$$

$$\leq \mathbb{P}\left(\sup_{0\leq t\leq n} |A_j(t) - \alpha_j t| + 1 + \alpha_j\tau' \geq \tau | \hat{a}_j(0) = 0\right)$$

$$\leq \mathbb{P}\left(\sup_{0\leq t\leq n} |A_j(t) - \alpha_j t| \geq (\tau/2) - 1 | \hat{a}_j(0) = 0\right),$$

where in the last inequality we have used the fact that $\alpha_j\tau' \leq \tau/2$. Since the bound above is uniform over $\tau'$ in the range $[0, \tau/(2\alpha_j)]$, it follows that

$$\mathbb{P}\left(\sup_{0\leq t\leq n} |A_j(t) - \alpha_j t| \geq \tau | a_j(0) \leq \tau/2, \hat{a}_j(0) = a_j\right)$$

$$\leq \mathbb{P}\bigg(\sup_{0\leq t\leq n} |A_j(t) - \alpha_j t| \geq (\tau/2) - 1 | \hat{a}_j(0) = 0\bigg).$$

Now, recall that $\tau = (\theta^{-1} n^{1/2} \log u)/3$, thus, plugging into the above bound, we get

$$\mathbb{P}\bigg(\sup_{0\leq t\leq n} |A_j(t) - \alpha_j t| \geq (\theta^{-1} n^{1/2} \log u)/6 - 1 | \hat{a}_j(0) = 0\bigg)$$

$$\leq \mathbb{P}\bigg(\sup_{0\leq t\leq n} |A_j(t) - \alpha_j t - (\Gamma_j^0)^{1/2} W_j^0(t)| \geq (\theta^{-1} n^{1/2} \log u)/12 - 1\bigg)$$

(51) $$+ \mathbb{P}\bigg(n^{-1/2} \sup_{0\leq t\leq n} |(\Gamma_j^0)^{1/2} W_j^0(t)| \geq (\theta^{-1} \log u)/12\bigg)$$

$$\leq n^{C_j} \exp(-(c_j/\theta)\sqrt{n} \log u + 1) + C_j' \exp(-c_j' \theta^{-2} (\log u)^2),$$

(52) $$= e n^{C_j} u^{-(c_j/\theta)\sqrt{n}} + C_j' \exp(-c_j' \theta^{-2} (\log u)^2),$$

for some finite positive constants $C_j, C_j'$ and $c_j, c_j'$, where the last step follows from the first displayed equation in Lemma A.2 and standard tail bounds on the maximum of a Brownian motion (see, e.g., [8], Lemma 1.6.1). For the second term in the RHS of (52), we have

$$\int_2^\infty C_j' e^{-c_j' \theta^{-2} (\log u)^2} \, du < \infty.$$

Let $n_0$ be the smallest integer such that $(c_j/\theta)\sqrt{n_0} > 1$. Then, for the first term on the RHS of (52), we have

$$\limsup_{n\to\infty} \int_2^\infty n^{C_j} u^{-(c_j/\theta)\sqrt{n}} \, du \leq \sup_{n\geq n_0} \int_2^\infty n^{C_j} u^{-(c_j/\theta)\sqrt{n}} \, du$$

$$= \sup_{n\geq n_0} \frac{n^{C_j}}{((c_j/\theta)\sqrt{n} - 1) 2^{(c_j/\theta)\sqrt{n}-1}} < \infty.$$

By the finiteness of the two bounds above, we obtain

$$\limsup_{n\to\infty} \int_2^\infty \mathbb{P}\bigg(\sup_{0\leq t\leq n} |A_j(t) - \alpha_j t| \geq (\theta^{-1} n^{1/2} \log u)/6 - 1$$

(53) $$|a_j(0) \leq \tau/(2\alpha_j), \hat{a}_j(0) = a_j\bigg) du$$

$$< \infty.$$

Putting this together with (49) and (50), we have

$$\limsup_{n\to\infty} \sup_{(z,a,v)\in\mathcal{X}} \mathbb{E}\bigg[\sup_{0\leq t\leq n} \exp(\theta n^{-1/2} |A_j(t) - \alpha_j t|) | \bar{Q}(0) = (z,a,v)\bigg]$$

$$\leq 2 + \limsup_n \sup_{(z,a,v)\in\mathcal{X}} \int_2^\infty \mathbb{P}\bigg(\sup_{0\leq t\leq n} |A_j(t) - \alpha_j t| \geq (\theta^{-1}n^{1/2}\log u)/3$$

(54)
$$|\bar{Q}(0) = (z,a,v)\bigg) du$$

$$< \infty.$$

Repeating these arguments for the service completion processes, we obtain similar bounds for the second term on the RHS of (47). To obtain the bound on the third term,

$$\mathbb{P}\bigg(\sup_{0\leq t\leq n}|R_j^i(S_i(t)) - p_{ij}S_i(t)| \geq (\theta^{-1}n^{1/2}\log u)/(3J)|\bar{Q}^n(0) = (z,a,v)\bigg),$$

we set again $\tau = \theta^{-1}n^{1/2}\log u/3$ and note that this term is bounded by

(55)
$$\mathbb{P}\bigg(\sup_{0\leq t\leq n+\tau/(J+1)}|R_j^i(\lfloor t\rfloor) - p_{ij}t| \geq \tau/J|\bar{Q}^n(0) = (z,a,v)\bigg)$$
$$+ \mathbb{P}\bigg(\sup_{0\leq t\leq n}|S_i(t) - \mu_i t| \geq \tau/(J+1)|\bar{Q}^n(0) = (z,a,v)\bigg).$$

We have shown above that, for the second probability in the sum above, the corresponding event $\sup_{0\leq t\leq n}|S_i(t) - \mu_i t| \geq \tau/(J+1)$ satisfies the bound of the form (54) with $S_i$ and $\mu_i$ replacing $A_j$ and $\alpha_j$, respectively. We now bound the first term in the sum above. Note that the corresponding event is independent from the event $\bar{Q}^n(0) = (z,a,v)$ since it only involves the routing process. We apply Lemma A.2 and, in particular, consider the corresponding coupled Brownian motion $W(t)$ with variance term $\Gamma_{i,j} = p_{ij}(1-p_{ij})$. We have

$$\mathbb{P}\bigg(\sup_{0\leq t\leq n+\tau/(J+1)}|R_j^i(\lfloor t\rfloor) - p_{ij}t| \geq \tau/J\bigg)$$
$$\leq \mathbb{P}\bigg(\sup_{0\leq t\leq n+\tau/(J+1)}|R_j^i(\lfloor t\rfloor) - p_{ij}t - (\Gamma_{i,j})^{1/2}W(n+\tau/(J+1))| \geq \tau/J\bigg)$$
$$+ \mathbb{P}\bigg(\sup_{0\leq t\leq n+\tau/(J+1)}|(\Gamma_{i,j})^{1/2}W(n+\tau/(J+1))| \geq \tau/J\bigg).$$

To bound the first term, we use the last part of Lemma A.2 and rewrite the term as

$$\mathbb{P}\bigg(\sup_{0\leq t\leq n+\tau/(J+1)}|R_j^i(\lfloor t\rfloor) - p_{ij}t - (\Gamma_{i,j})^{1/2}W(n+\tau/(J+1))|$$
$$\geq C_{i,j}^1\log(n+\tau/(J+1)) + (\tau/J - C_{i,j}^1\log(n+\tau/(J+1)))\bigg)$$

$$\leq C_{i,j}^2 \exp(-c_{i,j}\tau/J + c_{i,j}C_{i,j}^1 \log(n + \tau/(J+1)))$$
$$= (n + c_1' n^{1/2} \log u)^{c_2'} \exp(-c_3' n^{1/2} \log u),$$

where the last term is obtained by expanding $\tau$ as $\theta^{-1}n^{1/2}\log u/3$ and introducing the necessary constants $c_k', k = 1, 2, 3$. For all sufficiently large $n$, the last term is integrable with respect to $u$ and with the boundaries $\int_2^\infty$ uniformly in $n$. Namely,

$$\limsup_{n\to\infty} \int_2^\infty (n + c_1' n^{1/2} \log u)^{c_2'} \exp(-c_3' n^{1/2} \log u)\, du < \infty.$$

To bound the second probability in (55), we again use standard tail bounds on the maximum of a Brownian motion as in (51). The derivation of the bound is very similar, the only difference being that the upper limit $n$ in the supremum expression is replaced by $n + \tau/(J+1) = n + \theta^{-1}n^{1/2}\log u/(3(J+1))$. An easy check shows that the expression is still uniformly integrable in $u$ and the bound similar to (54) is obtained. We conclude that

$$\limsup_{n\to\infty} \sup_{(z,a,v)\in\mathcal{X}} \int_2^\infty \mathbb{P}\bigg( \sup_{0\leq t\leq n} |R_j^i(S_i(t)) - p_{ij}S_i(t)|$$
$$\geq (\theta^{-1}n^{1/2}\log u)/(3J)|\bar{Q}^n(0) = (z,a,v)\bigg) du < \infty.$$

To summarize, for any $0 \leq \theta \leq \theta^*$, we have

$$\limsup_{n\to\infty} \sup_{(z,a,v)\in\mathcal{X}} \mathbb{E}\bigg[\sup_{0\leq t\leq n} \exp(\theta n^{-1/2}|X_j^n(t) - x_j^n(t)|)|\bar{Q}^n(0) = (z,a,v)\bigg] < \infty.$$

Setting $\theta_1 < \theta^*/J$, we obtain

$$\limsup_{n\to\infty} \sup_{(z,a,v)\in\mathcal{X}} \mathbb{E}\bigg[\sup_{0\leq t\leq n} \exp(\theta n^{-1/2}\|X^n(t) - x^n(t)\|)|\bar{Q}^n(0) = (z,a,v)\bigg] < \infty,$$

which concludes the proof of (40). □

PROOF OF LEMMA A.2. The results of this lemma are immediate applications of functional strong approximations for random walks and the associated renewal processes, [7, 24, 25] (see also [11], Chapter 5). In particular, the following approximation results are established in the above referenced papers.

THEOREM 13. *Let $\xi_i, i = 1, 2, \ldots,$ be an i.i.d. sequence such that $\mathbb{E}[e^{\delta\xi_1}] < \infty$ for all sufficiently small $\delta > 0$. For every $t \in \Re_+$, define $X(t) = \sum_{i\leq \lfloor t\rfloor} \xi_i$ and the associated renewal process $Y(t) = \max\{k : \sum_{1\leq i\leq k} \xi_i < t\}$. There exist standard Brownian motions $W_1(t), W_2(t)$ defined on the same probability* space as $\xi_i$, and constants $C_1, C_2, C_3 > 0$ such that, for every $T > 1, u > 0$,

$$(56) \qquad \mathbb{P}\bigg(\sup_{0 \leq t \leq T} |X(t) - \nu^{-1}t - \sigma W(t)| \geq C_1 \log T + u\bigg) \leq C_2 e^{-C_3 u},$$

$$(57) \quad \mathbb{P}\bigg(\sup_{0 \leq t \leq T} |Y(t) - \nu t - (\nu)^{3/2} \sigma W(t)| \geq C_1 \log T + u\bigg) \leq C_2 e^{-C_3 u},$$

where $\nu = 1/\mathbb{E}[\xi_1]$ and $\sigma$ is the standard deviation of $\xi_1$.

The results of the lemma follow immediately from this theorem. The bounds corresponding to the arrival processes $A_j(t)$ and service processes $S_j(t)$ follow from (57) directly by setting $\nu = \alpha_j$ and $\nu = \mu_j$, respectively, and observing that $\nu^{3/2}\sigma = \alpha_j c_{a,j}^2$ for the $j$th arrival process and $= \mu_j c_{s,j}^2$ for the $j$th service process. The last bound corresponding to the routing process $R_j^k(\lfloor t \rfloor)$ follows from (56) by setting $\nu = p_{kj}$ and $\sigma = (\Gamma_{k,j})^{1/2}$. $\square$

**Acknowledgments.** The authors wish to thank Sean Meyn and Peter Glynn for many interesting and fruitful discussions. The authors also express their gratitude to the anonymous referees for their careful reading of the manuscript.

MIT SLOAN SCHOOL OF MANAGEMENT
CAMBRIDGE, MASSACHUSETTS 02139
USA
E-MAIL: gamarnik@watson.ibm.com

GRADUATE SCHOOL OF BUSINESS
COLUMBIA UNIVERSITY
NEW YORK, NEW YORK 10027
USA
E-MAIL: assaf@gsb.columbia.edu